\documentclass{amsart}
\usepackage{amscd,amssymb,graphicx}
\usepackage{graphicx}           
\author[Fialowski]{Alice Fialowski}
\address{
Alice Fialowski\\
E\"otv\"os Lor\'and University\\
Budapest, Hungary} \email{fialowsk@cs.elte.hu}
\author[Penkava]{Michael Penkava}
\address{
Michael Penkava\\
University of Wisconsin-Eau Claire\\
Eau Claire, WI 54702-4004} \email{penkavmr@uwec.edu}
\thanks{Research of the first author was partially supported by  OTKA grant
 K77757 and the Humboldt Foundation, the second author by grants from the
University of Wisconsin-Eau Claire.}

\newtheorem{thm}{Theorem}[section]

\theoremstyle{definition}

\def \ph{\varphi}

\def\GL{\mbox{\bf GL}}

\def \diag{\operatorname {diag}}

\def \ra{\rightarrow}

\def \hom{\mbox{\rm Hom}}
\def \ie{\hbox{\it i.e.}}

\def \gl{\mbox{$\mathfrak{gl}$}}
\def \tns{\otimes}

\def \mplus{+\cdots+}
\def \mcom{,\cdots,}
\def \k{\mbox{$\mathbb K$}}

\def \C{\mbox{$\mathbb C$}}
\def \Z{\mbox{$\mathbb Z$}}

\def\zt{\mbox{$\Z_2$}}

\def\ad{\operatorname{ad}}

\def\inv{^{-1}}
\def\d{d}

\def\im{\operatorname{Im}}

\def\L{L}

\def\m{\mbox{$\mathfrak m$}}

\def\coder{\operatorname{Coder}}

\def\and{\mbox{ \rm and }}
\def\T{\mathcal T}
\def\TV{\T(V)}
\def\TW{\mbox{$\T(W)$}}

\def\s#1{(-1)^{#1}}
\DeclareMathOperator*{\invlim}{\underleftarrow{\rm lim}}

\def\pha#1#2{\ph^{#1}_{#2}}

\def\psa#1#2{\psi^{#1}_{#2}}

\def\inv{^{-1}}

\def\dinf{\mbox{$d^\text{inf}$}}
\def\P{\mathbb P}

\input{xy}
\xyoption{all}

\begin{document}
\setlength{\multlinegap}{0pt}
\title[Non-nilpotent $4$-dimensional algebras]
{The Moduli space of $4$-dimensional non-nilpotent complex associative algebras}%

\address{}%
\email{}%

\thanks{}%

\begin{abstract}
In this paper, we study the moduli space of $4$-dimensional complex
associative algebras. We use extensions to compute the moduli space,
and then give a decomposition of this moduli space into strata consisting
of complex projective
orbifolds, glued together through jump deformations. Because the space of 4-dimensional
algebras is large, we only classify the non-nilpotent algebras in this paper.
\end{abstract}
\maketitle

\section{Introduction}

The classification of associative algebras was instituted by Benjamin Peirce
in the 1870's \cite{pie}, who gave a partial classification of the complex
associative algebras of dimension up to 6, although in some sense, one can
deduce the complete classification from his results, with some additional
work. The classification method relied on the following remarkable fact:
\begin{thm}
Every finite dimensional algebra which is not nilpotent contains a nontrivial
idempotent element.
\end{thm}
A nilpotent algebra $A$ is one which satisfies $A^n=0$ for some $n$, while an
idempotent element $a$ satisfies $a^2=a$.  This observation of Peirce eventually
leads to two important theorems in the classification of finite dimensional associative
algebras.  Recall that an algebra is said to be simple if it has no nontrivial
proper ideals, and it is not the trivial 1-dimensional nilpotent algebra over \k\, which is given
by the trivial product.
\begin{thm}[Fundamental Theorem of Finite Dimensional Associative Algebras]
Suppose that $A$ is a finite dimensional algebra over a field \k. Then $A$ has
a maximal nilpotent ideal $N$, called its radical.  If $A$ is not nilpotent, then $A/N$
is a semisimple
algebra, that is, a direct sum of simple algebras.
\end{thm}
In fact, in the literature, the definition of a semisimple algebra is often given as
one whose radical is trivial, and then it is a theorem that semisimple algebras are
direct sums of simple algebras. Moreover, when $A/N$ satisfies a property called separability
over \k, then $A$ is a semidirect product of its radical and a semisimple algebra.
Over the complex numbers, every semisimple algebra is separable. To apply this theorem
to construct algebras by extension, one uses the following characterization of simple algebras.
\begin{thm}
[Wedderburn] If $A$ is a finite dimensional algebra over \k, then $A$ is simple iff
$A$ is isomorphic
to a tensor product $M\tns D$, where $M=\mathfrak{gl}(n,\k)$ and
$D$ is a division algebra over \k.
\end{thm}
One can also say that $A$ is a matrix algebra with coefficients in a division algebra over \k.
An associative division algebra is a
unital associative algebra where every nonzero element has a multiplicative
inverse. (One has to modify this definition in the case of graded algebras, but we will
not address this issue in this paper.) Over the complex numbers, the only division algebra
is $\C$ itself, so Wedderburn's theorem says that the only simple algebras are the
matrix algebras. In particular, there is exactly one simple 4-dimensional complex associative
algebra, $\mathfrak{gl}(2,\C)$, while  there is one additional semisimple algebra,
the direct sum of 4 copies of $\C$.

According to our investigations, there are two basic prior approaches
to the classification. The first is the old paper by Peirce \cite{pie}
which attempts to classify all the nilpotent algebras, including
nonassociative ones. There are some evident mistakes in that paper, for
example, it gives a classification of the commutative nilpotent
associative algebras which contains nonassociative algebras as well.
The second approach \cite{mas} classifies the unital algebras only. It turns out
that classification of unital algebras is not sufficient.

Let us consider the unital algebra of one higher dimension which is
obtained by adjoining a multiplicative identity as the unital
enlargement of the algebra. Two nonisomorphic non-nilpotent algebras can
have isomorphic unital enlargements, so they cannot be recovered so
easily. Nevertheless, let us suppose that there were some efficient
method of constructing all unital algebras of arbitrary dimension, and
to determine their maximal nilpotent ideals. In that case, we could
recover all nilpotent algebras of dimension $n$ from their
enlargements. Moreover, to recover all algebras of dimension $n$, one
would only have to consider extensions of nilpotent algebras of
dimension $k$ by semisimple algebras of dimension $n-k$, where $0\leq k
\leq n$. Our method turns out to be efficient in
constructing extensions of nilpotent algebras by semisimple ones.

Thus, even if the construction of unital algebras could be carried out
simply, which is by no means obvious from the literature, one would
still need our methodology to construct most of the algebras. So the
role of our paper is to explore the construction method which leads to
the description of all algebras.

The main goal of this paper is to give a complete description of the
 moduli space of nonnilpotent 4-dimensional associative algebras,
including a computation of the miniversal deformation of every element.
We get the description with the help of extensions, which is the
novelty of our approach. The nilpotent cases will be classified in
another paper.
 We also give a canonical stratification
of the moduli space into projective orbifolds of a very simple type,
so that the strata are connected only by deformations
factoring through jump deformations, and the elements of a particular
stratum are given by neighborhoods determined by smooth
deformations.

The authors thank the referees for their useful comments.

\section{Construction of algebras by extensions}

In \cite{fp11}, the theory of extensions of an algebra $W$ by an algebra $M$ is
described. Consider the exact sequence
$$
0\ra M\ra V\ra W\ra 0
$$
of associative \k-algebras, so that $V=M\oplus W$ as a \k-vector space, $M$ is an
ideal in the algebra $V$, and $W=V/M$ is the quotient algebra. Suppose that
$\delta\in C^2(W)$ and $\mu\in C^2(M)$ represent the algebra structures on
$W$ and $M$ respectively. We can view $\mu$ and $\delta$ as elements of $C^2(V)$.
Let $T^{k,l}$ be the subspace of $T^{k+l}(V)$ given recursively
by
\begin{align*}
T^{0,0}&=\k\\
T^{k,l}&=M\tns T^{k-1,l}\oplus V\tns T^{k,l-1}
\end{align*}
Let
$C^{k,l}=\hom(T^{k,l},M)\subseteq C^{k+l}(V)$.
If we denote the algebra structure on $V$ by $d$, we have
$$
d=\delta+\mu+\lambda+\psi,
$$
where $\lambda\in C^{1,1}$ and $\psi\in C^{0,2}$. Note that in this notation,
$\mu\in C^{2,0}$. Then the condition that $d$ is associative:  $[d,d]=0$ gives the
following relations:
\begin{align*}
[\delta,\lambda]+\tfrac 12[\lambda,\lambda]+[\mu,\psi]&=0,
\quad\text{The Maurer-Cartan equation}\\
[\mu,\lambda]&=0,\quad\text{The compatibility condition}\\
[\delta+\lambda,\psi]&=0,\quad\text{The cocycle condition}
\end{align*}
Since $\mu$ is an algebra structure, $[\mu,\mu]=0$. Then if we define
$D_\mu$ by $D_\mu(\ph)=[\mu,\ph]$, then $D^2_\mu=0$.
Thus $D_\mu$ is a differential on $C(V)$.
Moreover $D_\mu:C^{k,l}\ra C^{k+1,l}$. Let
\begin{align*}
Z_\mu^{k,l}&=\ker(D_\mu:C^{k,l}\ra C^{k+1,l}),\quad\text{the $(k,l)$-cocycles}\\
B_\mu^{k,l}&=\im(D_\mu:C^{k-1,l}\ra C^{k,l}),\quad\text{the $(k,l)$-coboundaries}\\
H_\mu^{k,l}&=Z_\mu^{k,l}/B_\mu^{k,l},\quad\text{the $D_u$ $(k,l)$-cohomology}
\end{align*}

Then the compatibility condition means that $\lambda\in Z^{1,1}$.
If we define $D_{\delta+\lambda}(\ph)=[\delta+\lambda,\ph]$, then it is not
true that $D^2_{\delta+\lambda}=0$, but
$D_{\delta+\lambda}D_\mu=-D_{\mu}D_{\delta+\lambda}$, so that $D_{\delta+\lambda}$ descends
to a map $D_{\delta+\lambda}:H^{k,l}_\mu\ra H^{k,l+1}_\mu$, whose square is zero, giving
rise to the $D_{\delta+\lambda}$-cohomology $H^{k,l}_{\mu,\delta+\lambda}$.
Let the pair $(\lambda,\psi)$ give rise to a codifferential $d$, and $(\lambda,\psi')$
give rise to another codifferential $d'$. Then if we express $\psi'=\psi+\tau$, it is
easy to see that $[\mu,\tau]=0$, and $[\delta+\lambda,\tau]=0$, so that the image $\bar\tau$
of $\tau$ in $H^{0,2}_\mu$ is a $D_{\delta+\lambda}$-cocycle, and thus $\tau$ determines
an element $\{\bar\tau\}\in H^{0,2}_{\mu,\delta+\lambda}$.

If $\beta\in C^{0,1}$, then $g=\exp(\beta):\TV\ra\TV$ is given by
$g(m,w)=(m+\beta(w),w)$. Furthermore $g^*=\exp(-\ad_{\beta}):C(V)\ra C(V)$ satisfies
$g^*(d)=d'$, where $d'=\delta+\mu+\lambda'+\psi'$ with
$\lambda'=\lambda+[\mu,\beta]$ and $\psi'=\psi+[\delta+\lambda+\tfrac12[\mu,\beta],\beta]$.
In this case, we say that $d$ and $d'$ are equivalent extensions in the restricted sense.
Such equivalent extensions are also equivalent as codifferentials on $\TV$.
Note that $\lambda$
and $\lambda'$ differ by a $D_\mu$-coboundary, so $\bar\lambda=\bar\lambda'$ in
$H^{1,1}_\mu$. If $\lambda$ satisfies the MC equation for some $\psi$, then
any element $\lambda'$ in $\bar\lambda$ also gives a solution of the MC equation,
for the $\psi'$ given above. The cohomology classes of those $\lambda$ for which
a solution of the MC equation exists determine distinct restricted equivalence classes
of extensions.

Let $G_{M,W}=\GL(M)\times\GL(W)\subseteq\GL(V)$. If $g\in G_{M,W}$ then $g^*:C^{k,l}\ra
C^{k,l}$, and $g^*:C^k(W)\ra C^k(W)$, so $\delta'=g^*(\delta)$ and $\mu'=g^*(\mu)$
are codifferentials on $\T(M)$ and $\TW$ respectively.
The group $G_{\delta,\mu}$ is the
subgroup of $G_{M,W}$ consisting of those elements $g$ such that $g^*(\delta)=\delta$
and $g^*(\mu)=\mu$. Then $G_{\delta,\mu}$ acts on
the restricted equivalence classes of extensions, giving the equivalence classes
of general extensions. Also $G_{\delta,\mu}$
acts on $H^{k,l}_\mu$, and induces an action on the classes $\bar\lambda$ of $\lambda$
giving a solution to the MC equation.

Next, consider the group $G_{\delta,\mu,\lambda}$ consisting
of the automorphisms $h$ of $V$ of the form $h=g\exp(\beta)$, where
$g\in G_{\delta,\mu}$, $\beta\in C^{0,1}$ and $\lambda=g^*(\lambda)+[\mu,\beta]$.
If $d=\delta+\mu+\lambda+\psi+\tau$, then $h^*(d)=\delta+\mu+\lambda+\psi+\tau'$ where
\begin{equation*}
\tau'=g^*(\psi)-\psi+[\delta+\lambda-\tfrac12[\mu,\beta],\beta]+g^*(\tau).
\end{equation*}
Thus the group $G_{\delta,\mu,\lambda}$ induces an action on $H^{0,2}_{\mu,\delta+\lambda}$
given by $\{\bar\tau\}\ra\{\overline{\tau'}\}$.

The general group of equivalences of extensions of the algebra structure $\delta$ on $W$
by the algebra structure $\mu$ on $M$ is given by the group of automorphisms of $V$ of
the form $h=\exp(\beta)g$, where $\beta\in C^{0,1}$ and $g\in G_{\delta,\mu}$. Then we
have the following classification of such extensions up to equivalence.
\begin{thm}[\cite{fp11}]
The equivalence classes of extensions of $\delta$ on $W$ by $\mu$ on $M$ is classified
by the following:
\begin{enumerate}
\item Equivalence classes of $\bar\lambda\in H^{1,1}_\mu$ which satisfy the MC equation
\begin{equation*}
[\delta,\lambda]+\tfrac12[\lambda,\lambda]+[\mu,\psi]=0
\end{equation*}
for some $\psi\in C^{0,2}$, under the action of the group $G_{\delta,\mu}$.
\item Equivalence classes of $\{\bar\tau\}\in H^{0,2}_{\mu,\delta+\lambda}$ under the
action of the group $G_{\delta,\mu,\lambda}$.
\end{enumerate}
\end{thm}
Equivalent extensions will give rise to equivalent algebras on $V$, but it may
happen that two algebras arising from nonequivalent extensions are equivalent.
This is because the group of equivalences of extensions is the group of invertible
block upper
triangular matrices on the space $V=M\oplus W$, whereas the the equivalence
classes of algebras on $V$ are given by the group of all invertible
matrices, which is larger.

The fundamental theorem of finite dimensional algebras allows us to restrict our
consideration of extensions to two cases. First, we can consider those extensions
where $\delta$ is a semisimple algebra structure on $W$, and $\mu$ is a nilpotent
algebra structure on $M$. In this case, because we are working over $\C$, we can
also assume that $\psi=\tau=0$. Thus the classification of the extension reduces
to considering equivalence classes of $\lambda$.

Secondly, we can consider extensions
of the trivial algebra structure $\delta=0$ on a 1-dimensional space $W$ by
a nilpotent algebra $\mu$. This
is because a nilpotent algebra has a codimension 1 ideal $M$, and the restriction
of the algebra structure to $M$ is nilpotent. However, in this case, we cannot assume
that $\psi$ or $\tau$ vanish,
so we need to use the classification theorem above to determine the
equivalence classes of extensions. In many cases, in solving the MC equation for
a particular $\lambda$, if there is any $\psi$ yielding a solution, then $\psi=0$
also gives a solution, so the action of $G_{\delta,\mu,\lambda}$ on $H^{0,2}_\mu$
takes on a simpler form than the general action we described above.

In addition to the complexity which arises because we cannot take the cocycle term
$\psi$ in the extension to be zero, there is another issue that complicates the construction
of the extensions. If an algebra is not nilpotent, then it has a maximal nilpotent ideal
which is unique, and it can be constructed as an extension of a semisimple algebra by this
unique ideal. Both the semisimple and nilpotent parts in this construction are completely
determined by the algebra. Therefore, a classification of extensions up to equivalence of
extensions will be sufficient to classify the algebras. This means that the equivalence classes
of the module structure $\lambda$ determine the algebras up to isomorphism.

For nilpotent algebras, we don't have this assurance. The same algebra structure may arise by
extensions of the trivial algebra structure on a 1-dimensional space by two different nilpotent
algebra structures on the same $n-1$-dimensional space.

In addition, the deformation theory of the nilpotent algebras is far more involved than the deformation
theory of the nonnilpotent algebras. Thus, we decided to discuss the nilpotent 4-dimensional complex
algebras in a separate paper. In this paper, we only look at extensions of semisimple algebras by
nilpotent algebras, which is precisely what is necessary to classify
all non-nilpotent algebras.

\section{Associative algebra structures on a 4-dimensional vector space}
Denote the basis elements of a 4-dimensional associative algebra by
$f_1,f_2,f_3, f_4$ and let $\psi_k^{ij}$ denote the product
$f_if_j=f_k$.
 We will
recall the classification of algebras on  a 2-dimensional space given in \cite{bdhoppsw2},
and the classification of algebras on a 3-dimensional space given in \cite{fpp1}.
\begin{table}[h]
\begin{center}
\begin{tabular}{lccccc}
Codifferential&$H^0$&$H^2$&$H^1$&$H^3$&$H^4$\\ \hline \\
$d_1=\psi_1^{11}+\psi_2^{22}$&2&0&0&0&0\\
$d_2=\psi_2^{22}+\psi^{12}_1$&0&0&0&0&0\\
$d_3=\psi_2^{22}+\psi_1^{21}$&0&0&0&0&0\\
$d_4=\psi_2^{22}+\psi_1^{12}+\psi_1^{21}$&2&1&1&1&1\\
$d_5=\psi_2^{22}$&2&1&1&1&1\\
$d_6=\psi_1^{22}$&2&2&2&2&2\\\\ \hline
\end{tabular}
\end{center}
\label{Table 1}
\caption{Two dimensional complex associative algebras and their
cohomology}
\end{table}

\begin{table}[h]
\begin{center}
\begin{tabular}{lccccc}
Codifferential&$H^0$&$H^2$&$H^1$&$H^3$&$H^4$\\ \hline \\
$d_1=\psa{33}3+\psa{22}2+\psa{11}1$&3&0&0&0&0\\
$d_{2}=\psi_2^{22}+\psi_3^{33}+\psi_1^{21}+\psi_1^{13}$&1&0&0&0&0\\
$d_{3}=\psi_2^{22}+\psi_3^{33}+\psi_1^{12}$&1&0&0&0&0\\
$d_{4}=\psi_2^{22}+\psi_3^{33}+\psi_1^{21}$&1&0&0&0&0\\
$d_{5}=\psi_2^{22}+\psi_3^{33}+\psi_1^{21}+\psi_1^{12}$&3&1&1&1&1\\
$d_{6}=\psi_2^{22}+\psi_3^{33}$&3&1&1&1&1\\
$d_{7}=\psi_3^{33}+\psi_1^{22}+\psi_1^{31}+\psi_2^{32}+\psi_1^{13}+\psi_2^{23}$&3&2&2&2&2\\
$d_{8}=\psi_3^{33}+\psi_1^{22}$&3&2&2&2&2\\
$d_{9}=\psi_3^{33}+\psi_1^{31}+\psi_2^{32}$&0&3&0&0&0\\
$d_{10}=\psi_3^{33}+\psi_1^{13}+\psi_2^{23}$&0&3&0&0&0\\
$d_{11}=\psi_3^{33}+\psi_1^{31}+\psi_2^{23}$&0&1&0&1&0\\
$d_{12}=\psi_3^{33}+\psi_1^{13}+\psi_2^{32}+\psi_2^{23}$&1&1&1&1&1\\
$d_{13}=\psi_3^{33}+\psi_1^{31}+\psi_2^{32}+\psi_2^{23}$&1&1&1&1&1\\
$d_{14}=\psi_3^{33}+\psi_2^{32}$&1&1&2&2&2\\
$d_{15}=\psi_3^{33}+\psi_2^{23}$&1&1&2&2&2\\
$d_{16}=\psi_3^{33}+\psi_2^{32}+\psi_2^{23}$&3&2&2&2&2\\
$d_{17}=\psi_3^{33}+\psi_1^{31}+\psi_1^{13}+\psi_2^{32}+\psi_2^{23}$&3&4&6&12&24\\
$d_{18}=\psi_3^{33}$&3&4&8&16&32\\
$d_{19}=\psi_2^{13}+\psi_2^{31}+\psi_1^{33}$&3&3&3&3&3\\
$d_{20}(0:0)=\psi_2^{33}$&3&5&9&17&33\\
$d_{20}(1:0)=\psi_2^{13}+\psi_2^{33}$&1&2&5&8&11\\
$d_{20}(1:1)=\psi_2^{13}+\psi_2^{31}+\psi_2^{33}$&3&4&5&7&8\\
$d_{20}(1:-1)=\psi_2^{13}-\psi_2^{31}+\psi_2^{33}$&1&2&3&4&5\\
$d_{20}(p:q)=\psi_2^{13}p+\psi_2^{31}q+\psi_2^{33}$&1&2&3&3&4\\
$d_{21}=\psi_2^{13}-\psi_2^{31}$&1&4&5&8&9\\\\ \hline
\end{tabular}
\end{center}
\label{coho3 table}
\caption{Three dimensional complex associative algebras and their
cohomology}
\end{table}



Actually, we only need to know the nilpotent algebras from lower dimensions as well
as the semisimple algebras. In dimension 1, there is one nontrivial algebra structure
$d_1=\psa{11}1$, which is just complex numbers $\C$.

Thus, in dimension 2, the algebra $d_1=\psi_1^{11}+\psi_2^{22}$ is semisimple, while
the algebra $d_6=\psi_1^{22}$ is nilpotent. These are the only algebras of dimension 2 (other than
the trivial algebra) which play a role in the construction of 4-dimensional algebras by
extensions. The algebra $d_1$ is just the direct sum $\C^2=\C\oplus\C$.

In the case of 3-dimensional algebras, only $d_1=\psi_1^{11}+\psi_2^{22}+\psi_3^{33}$ is semisimple, and only the
algebras $d_{19}=\psi_2^{13}+\psi_2^{31}+\psi_1^{33}$, $d_{20}(p:q)=\psi_2^{13}p+\psi_2^{31}q+\psi_2^{33}$, and $d_{21}=\psi_2^{13}-\psi_2^{31}$ are nilpotent. The algebra $d_1$ is
just the direct sum of three copies of $\C$.

Note that $d_{20}(p:q)$ is
a family of algebras parameterized by the projective orbifold $\P^1/\Sigma_2$. By this
we mean that the algebras $d_{20}(p:q)$ and $d_{20}(tp:tq)$ are isomorphic if $t\ne0$,
which gives the projective parameterization,
and that the algebras $d_{20}(p:q)$ and $d_{20}(q:p)$ are also isomorphic, which gives
the action of the group $\Sigma_2$ on $\P^1$.

In constructing the elements of the moduli space by extensions, we need to consider
three possibilities, extensions of the semisimple algebra structure on a 3-dimensional space
$W$ by the
trivial algebra structure on a 1-dimensional space $M$,
extensions of the semisimple algebra structure on a 2-dimensional space by a nilpotent
algebra on a 2-dimensional space, and extensions of either the simple or the trivial
1-dimensional algebras by a nilpotent 3-dimensional algebra.
\medskip

Consider the general setup,
where an $n$-dimensional space
$W=\langle f_{m+1},\dots f_{m+n}\rangle$ is extended
by an $m$-dimensional space $M=\langle f_1\mcom f_m\rangle$. Then the module structure is
of the form
$$
\lambda=\psa{kj}i(L_k)^i_j+\psa{jk}i(R_k)^i_j,\quad i,j=1,\dots m, k=m+1\dots m+n,
$$
and we can consider $L_k$ and $R_k$ to be $m\times m$ matrices. Then we can express
the bracket $\tfrac12[\lambda,\lambda]$, which appears in the MC equation in terms
of matrix multiplication.

\begin{align}
\tfrac12[\lambda,\lambda]=
\psa{jkl}i(R_lR_k)^i_j+\psa{kjl}i(L_kR_l-R_kL_l)^i_j-\psa{klj}i(L_kL_l)^i_j,
\label{lambda-lambda}
\end{align}
where $i,j=1,\dots m$, and $k,l=m+1,\dots m+n$.

Next, suppose that $\delta=\psa{m,m}m\mplus \psa{m+n,m+n}{m+n}$ is the semisimple
algebra structure $\C^n$ on $W$. Then we can also express $[\delta,\lambda]$ in
terms of matrix multiplication.
\begin{equation}
[\delta,\lambda]=\psa{kkj}i(L_k)^i_j-\psa{jkk}i(R_k)^i_j.
\label{delta-lambda}
\end{equation}
Since $\delta$ is semisimple, one can ignore the cocycle $\psi$ in constructing an extension,
so the MC equation is completely determined by the equations (\ref{delta-lambda}) and
(\ref{lambda-lambda}), so we obtain the conditions.
Therefore, the MC equation holds precisely when
\begin{align*}
L_k^2=L_k,\quad R_k^2=R_k,\quad L_kL_l=R_kR_l=0\text{ if $k\ne l$},\quad L_kR_l=R_lL_k.
\end{align*}
As a consequence, both $L_k$ and
$R_k$ must be commuting nondefective matrices whose eigenvalues are either 0 or 1, which
limits the possibilities.
Moreover, it can be shown that $G_{\delta}$, the group of automorphisms of $W$ preserving
$\delta$ is just the group of permutation matrices. Thus if $G=\diag(G_1,G_2)$ is
a block diagonal element of $G_{\delta,\mu}$, the matrix $G_2$ is a permutation matrix.
The action of $G$ on $\lambda$ is given by simultaneous conjugation of the matrices
$L_k$ and $R_k$ by $G_1$, and a simultaneous permutation of the $k$-indices determined
by the permutation associated to $G_2$.

When $\mu$ is zero, this is  the entire story. When $\mu\ne0$, the matrices $G_1$ are
required to preserve $\mu$, and the compatibility condition $[\mu,\lambda]$ also complicates
the picture.

It is important to note that given an $m$ and a nilpotent element $\mu$ on an $m$-dimensional
space $M$, there is an $n$ beyond which the extensions of the semisimple codifferential on an
$N$ dimensional space with $N$ greater than $n$ are simply direct sums of the extensions
of the $n$-dimensional semisimple algebra $\C^n$ and the semisimple algebra $\C^{N-n}$.
We say that the extension theory becomes stable at $n$. Moreover, the deformation picture
stabilizes as well.

In higher dimensions, there are semisimple algebras which are not of the form $\C^n$. Also,
as $m$ increases, the complexity of the nontrivial nilpotent elements $\mu$ increases as
well. In dimension 4, there is a simple algebra, $\gl(2,\C)$, represented by the
codifferential $d_1$, and a semisimple algebra $\C^4$, represented by the algebra $d_2$.
All other 4-dimensional nonnilpotent algebras are extensions of a semisimple algebra of
the type $\C^n$, for $n=1,2,3$.

\section{Extensions of the 3-dimensional semisimple algebra $\C^3$
by the 1-dimensional trivial
algebra $\C_0$}
Let $W=\langle f_2, f_3,f_4\rangle$ and $M=\langle f_1\rangle$.
The matrices $L_k$ and $R_k$ determining  $\lambda$ are $1\times1$ matrices, in
other words, just numbers; in fact, they are either 0 or 1. By applying a permutation
to the indices $2,3,4$, we can assume that either all the $L_k$ vanish, or $L_2=1$
and both $L_3$ and $L_4$ vanish. In the first case, either $R_2=1$ or $R_2=0$ and
$R_3=R_4=0$. In the second case, we can either have $R_2=1$ and $R_3=0$, $R_2=0$ and
$R_3=1$, or both $R_2$ and $R_3$ vanish. In all three cases, $R_4=0$.
Note that in all of these solutions, we can assume that $L_4=R_4=0$.  For extensions
by a 1-dimensional space $M$, the extension picture stabilizes at $n=2$, and we are
looking at $n=3$. Thus the five solutions for $\lambda$ here, which
give  the codifferentials $d_3\mcom d_7$ correspond to the five 3-dimensional codifferentials
$d_2\mcom d_6$.

\section{Extensions of the 2-dimensional semisimple algebra $\C^2$
by a 2-dimensional nilpotent algebra}

Let $W=\langle f_3,f_4\rangle$ and $M=\langle f_1,f_2\rangle$.

There are two choices of
$\mu$ in this case, depending on whether the algebra structure on $M$ is  the
trivial or nontrivial nilpotent
structure. Although we cannot calculate $G_{\delta,\mu}$ without knowing $\mu$, we can say
that the matrix $G_2$ in the expression above for an element of $G_{M,W}$ must be one
of the two permutation matrices.
\subsection{Extensions by the nontrivial nilpotent algebra}
In this case $\mu=\psa{22}1$. In order for $[\mu,\lambda]=0$, using
equations (\ref{delta-lambda}) and (\ref{lambda-lambda}), we must have
$$ (L_k)^2_1=(R_k)^2_1=0, (L_k)^1_1=(L_k)^2_2=(R_k)^1_1=(R_k)^2_2,
$$
for all $k$. It follows that $L_k$ and $R_k$ are upper triangular matrices with
the same values on the diagonal, and since they are also nondefective matrices, they
must be diagonal, and therefore are either both equal to the identity or both the zero
matrix. It follows that by applying a permutation, we obtain either the solution
$\lambda=0$ or $L_3=R_4=I$, and $L_k=R_k=0$ for $k>3$. In fact, we see that this
case stabilizes when $n=1$, and we are looking at the case $n=2$. Thus the two solutions
$d_8$ and $d_9$ correspond to the three dimensional algebras $d_7$ and $d_8$. In fact,
$d_7$ arises by the following consideration. Given an algebra on an $n$-dimensional space,
there is an easy way to extend it to a unital algebra on an $n+1$-dimensional space, by taking any
vector not in the original space and making it play the role of the identity. This is how $d_7$ arises.
The algebra $d_8$ also arises in a natural way as the direct sum of the algebra structures $\delta$ and $\mu$.
For this $\mu$, these are the only such structures which arise, and this is somewhat typical.

\subsection{Extensions by the trivial nilpotent algebra}
In this case, $L_k$ and $R_k$ are $2\times2$ matrices.
The nontrivial permutation has the effect of
interchanging $L_3$ and $L_4$ as well as $R_3$ and $R_4$. The matrix $G_1$ acts on
all four of the matrices by simultaneously conjugating them.

By permuting if necessary, one can assume that $L_3$ is either a nonzero matrix,
or both $L_3$ and $L_4$ vanish. Moreover, by conjugation in case $L_3$ is not
the identity or the zero matrix, we have $L_3=\diag(1,0)$, which we will denote by
$T$. Now if $L_3=I$, then since $L_3L_4=0$, we must have $L_4=0$, but if $L_3=T$,
then the condition $L_3L_4=0$ forces to either be
$B=\diag(0,1)$ or 0. A similar analysis applies to the $R$ matrices.
Let us consider a case by case analysis.

If $L_3=I$ then $L_4=0$. Since $I$ is invariant
under conjugation, we can still apply a conjugation to put $R_3$ in the form
$I$, $T$ or $0$. If $R_3=I$, then $R_3=0$. If $R_3=T$, then either $R_4=B$ or
$R_4=0$. If $R_3=0$, then $R_4$ may equal $I$, $T$ or $0$. This gives six
solutions.

Next, assume $L_3=T$ and $L_4=B$.
Since we have used up the conjugation in putting $L_3$ and $L_4$ in
diagonal form, we can only use the fact that since $R_3$ and $R_4$
commute with $L_3$ and $L_4$, they can be simultaneously diagonalized, so we may
assume they are diagonal. Thus $R_3$ is either $I$, $T$ $B$ or $0$.
If $R_3=I$ then $R_4=0$. If $R_3=T$ then
$R_4=B$ or $R_4=0$. If $R_3=D$ then $R_4=T$ or $R_4=0$. If $R_3=0$, then $R_4$ is
either $I$, $T$, $B$ or $0$. This gives 9 possibilities, but there is one more thing
which we have to be careful of. A certain conjugation interchanges $T$ and $D$, so that if we
first apply the nontrivial permutation and then the conjugation which interchanges $T$ and $B$,
we find that the $L_3=T$, $L_4=B$, $R_3=I$ and $R_4=0$ is the same as if $R_3=0$ and
$R_4=I$. Similarly, $R_3=T$ and $R_4=0$ transforms to $R_3=0$ and $R_4=B$. Finally,
$R_3=B$ and $R_4=0$ transforms to $R_3=0$ and $R_4=B$. Thus instead of 9 cases, we
only obtain 6.

If $L_3=T$ and $L_4=0$, then we obtain the same 9 cases for $R_3$ and $R_4$ as when
$L_3=T$ and $L_4=B$, except this time, there are no hidden symmetries, so we get
exactly 9 cases.

Finally, when $L_3=L_4=0$, if $R_3=I$ then $R_4=0$, while if $R_3=T$, then $R_4=B$
or $R_4=0$, and if $R_3=0$, then $R_4=0$, giving 4 more cases.

This gives a total of 25 nonequivalent extensions, and they are also all nonequivalent as
algebras, corresponding to $d_{10},\dots d_{34}$. In this case, $n=2$ is not the
stable case. It is not hard to see that $n=4$ gives the stable case, corresponding to
the 6-dimensional moduli space.
\section{Extensions of the 1-dimensional simple algebra $\C$
by a 3-dimensional nilpotent algebra}
Here $M=\langle f_1,f_2,f_3\rangle$ and $W=\langle f_4\rangle$, and  $L_4$
and $R_4$ are $3\times 3$ matrices, which for simplicity, we will just denote by $L$ and $R$.
Elements in $C^{0,1}$ are of the form Let $\beta=\pha{4}1b_1+\pha{4}2b_2+\pha43b_3$.
\subsection{Extensions by the nilpotent algebra $\mu=\psa{13}2+\psa{31}2+\psa{33}1$}
In order for $[\mu,\lambda]=0$, using equations (\ref{delta-lambda})
and (\ref{lambda-lambda}), we must have
$$L=\left[\begin{matrix}L^1_1&0&L^1_3\\L^1_3&L^1_1&L^2_3\\0&0&L^1_1\end{matrix}
\right],
R=\left[\begin{matrix}L^1_1&0&L^1_3\\L^1_3&L^1_1&R^2_3\\0&0&L^1_1\end{matrix}
\right].$$
Taking into account the MC equation, we obtain that $L$ and $R$ must be diagonal matrices,
so we only get two solutions, depending on whether $L$ and $R$ both vanish or are both
equal to the identity matrix. Note that this is the stable case. It corresponds to the
codifferentials $d_{35}$ and $d_{36}$. Notice that we obtain one unital algebra $d_{35}$ and one algebra
which is a direct sum, $d_{36}$.
\subsection{Extensions by the nilpotent algebra $\mu=\psa{13}2p+\psa{31}2q+\psa{33}2$}
Here the situation depends on the projective coordinate $(p:q)$, which is parameterized by $\P^1/\Sigma_2$.
There are special cases when $p=q=0$ or $p=1$ and $q=0$. These three cases arise
from the compatibility condition $[\mu,\lambda]=0$, which generically has one solution,
but has additional solutions when either $p$ or $q$ vanishes, and when both $p$ and $q$ vanish.
\subsubsection{The generic case}
In this case, in order for $[\mu,\lambda]=0$, we must have
$$L=\left[\begin{matrix}L^1_1&0&0\\L^2_1&L^1_1&L^2_3\\0&0&L^1_1\end{matrix}
\right],
R=\left[\begin{matrix}L^1_1&0&0\\R^2_1&L^1_1&R^2_3\\0&0&L^1_1\end{matrix}
\right].$$

Taking into account the MC equation, we obtain that $L$ and $R$ must be diagonal matrices,
so we only get two solutions, depending on whether $L$ and $R$ both vanish or are both
equal to the identity matrix. Note that this is the stable case. It corresponds to the
algebras $d_{37}(p:q)$ and $d_{38}(p:q)$. Note that both of these are families
parameterized by $\P^1/\Sigma_2$.
\subsubsection{The case $p=1$, $q=0$}
In this case, in order for $[\mu,\lambda]=0$, we must have
\begin{equation*}
L=\left[ \begin {array}{ccc}
L^1_1&0&L^1_1-L^3_3\\\noalign{\medskip}
L^2_1&L^1_1&L^2_3\\\noalign{\medskip}
0&0&L^3_3
\end {array} \right] , R=\left[\begin {array}{ccc}
L^3_3&0&L^3_3-R^3_3\\\noalign{\medskip}
R^2_1&R^3_3&R^2_3\\\noalign{\medskip}
0&0&R^3_3
\end {array} \right].
\end{equation*}
Since $[\mu,\beta]=\psa{43}2(b_1+b_3)+\psa{14}2b_3+\psa{34}2b_3$, we can further
assume that $L^2_3=0$. Moreover, the eigenvalues of $L$ are $L^1_1$ and $L^3_3$,
while those of $R$ are $L^3_3$ and $R^3_3$, and these numbers must be either 0
 or 1,
yielding 8 possibilities. In fact, each one of these 8 choices corresponds to a
 solution
solution of the MC equation. Let
\begin{equation*}
T_1=\left[\begin{matrix}
1&0&1\\
0&1&0\\
0&0&0
\end{matrix}\right],
T_2=\left[\begin{matrix}
0&0&-1\\
0&0&0\\
0&0&1
\end{matrix}\right],
B_1=\left[\begin{matrix}
1&0&1\\
0&0&0\\
0&0&0
\end{matrix}\right],
B_2=\left[\begin{matrix}
0&0&-1\\
0&1&0\\
0&0&1
\end{matrix}\right].
\end{equation*}
Then the 8 solutions are $L=I,R=I$;\enspace $L=0,R=0$;\enspace
$L=I,R=B_1$;\enspace
$L=T_1,R=B_2$;\enspace
$L=T_2,R=B_2$;\enspace $L=0,R=B_2$;\enspace $L=T_2,R=B_1$ and $L=T_1,R=0$, corresponding to
the codifferentials
$d_{37(1:0)}$,\enspace $d_{38(1: 0)}$,\enspace $d_{39}\mcom d_{44}$.

This is not the stable case, because given an $L,R$ pair above,
there is another such pair, which satisfies the requirements that the products of the
$L$ matrices vanish, the products of the $R$ matrices vanish, and the $L$ and $R$ matrices commute.
In fact, it is not hard to see that $n=2$ gives the stable case, which will occur for
5-dimensional algebras.
\subsubsection{The case $p=0$, $q=0$}
In this case, in order for $[\mu,\lambda]=0$, we must have
\begin{equation*}
L=\left[ \begin {array}{ccc}
L^1_1&0&L^1_3\\\noalign{\medskip}
L^2_1&L^3_3&L^2_3\\\noalign{\medskip}
0&0&L^3_3
\end {array} \right] , R=\left[\begin {array}{ccc}
R^1_1&0&R^1_3\\\noalign{\medskip}
R^2_1&L^3_3&R^2_3\\\noalign{\medskip}
0&0&L^3_3
\end {array} \right].
\end{equation*}

Since $[\mu,\beta]=\psa{43}2b_3+\psa{34}2b_3$, we can further
assume that $L^2_3=0$. Moreover, the eigenvalues of $L_4$ are $L^1_1$ and $L^3_3$,
while those of $R_4$ are $L^3_3$ and $R^3_3$, and these numbers must be either 0 or 1,
yielding 8 possibilities. In fact, each one of these 8 choices corresponds to a unique
solution of the MC equation, which has some parameters. However, at this point we still
have not taken into account the action of $G_{\delta,\mu}$. An element in
$G_\mu$ is a matrix of the form $G=
\left[ \begin {smallmatrix}
g^1_1&0&g^1_3\\\noalign{\medskip}
g^2_1&(g^3_3)^2&g^2_3\\\noalign{\medskip}
0&0&g^3_3
\end {smallmatrix} \right]$,
and it acts on $\lambda$ by conjugating $L$ and $R$ simultaneously. This action is sufficient
to eliminate the parameters in the solutions for $L$ and $R$. Let $T=\diag(1,0,0)$ and
$B=\diag(0,1,1)$.

The 8 solutions are $L=I,R=I$, $L=0,R=0$, $L=I,R=B$, $L=T,R=T$,
$L=B,R=I$, $L=0,R=T$, $L=B,R=B$ and $L=T,R=0$, corresponding to
the codifferentials $d_{37(0: 0)}$, $d_{38(0: 0)}$,
$d_{45}\mcom d_{50}$.

This is not the stable case, because given an $L,R$ pair above,
there is another such pair, which satisfies the requirements that the products of the
$L$-s vanish, the products of the $R$-s vanish, and the $L$ and $R$ matrices commute.
In fact, it is not hard to see that $n=2$ gives the stable case, which will occur for
5-dimensional algebras.
\subsection{Extensions by the nilpotent algebra $\mu=\psa{13}2-\psa{31}2$}
In order for $[\mu,\lambda]=0$, we must have
$$L=\left[ \begin {matrix}
L^3_3&0&0\\\noalign{\medskip}
L^2_1&L^3_3&L^2_3\\\noalign{\medskip}
0&0&L^3_3
\end {matrix} \right] , \left[ \begin {matrix}
L^3_3&0&0\\\noalign{\medskip}
R^2_1&L^3_3&R^2_3\\\noalign{\medskip}
0&0&L^3_3
\end {matrix} \right]$$
Taking into account the MC equation, we obtain that $L$ and $R$ must be diagonal matrices,
so we only get two solutions, depending on whether $L$ and $R$ both vanish or are both
equal to the identity matrix. Note that this is the stable case. It corresponds to the
codifferentials $d_{51}$, which is the unital extension, and $d_{52}$, which is the direct sum extension.
\subsection{Extensions by the trivial nilpotent algebra}
Since $\mu=0$, we don't get any restrictions on $\lambda$ from the compatibility
condition, but since $G_{\mu}=\GL(3,\C)$, we can assume that $L$ is in Jordan normal
form, and since $L$ is nondefective, this implies that $L$ is diagonal. From
this it follows that $L$ can only be one of $I$, $T_1=\diag(1,1,0)$ $T_2=\diag(1,0,0)$
or 0.

When $L=I$ or $L=0$, it is invariant under conjugation, so we may conjugate $R$ to be one
of the same 4 matrices $I$, $T_1$, $T_2$ or 0.  When $L=T_1$, $R$ can also be conjugated
to make it diagonal, and we obtain  that $R$ is one of the six matrices $I$, $0$, $T_1$, $T_2$,
$B_1=\diag(1,0,1)$ or $B_2=\diag(0,0,1)$. When $L=T_2$, $R$ can again be conjugated to make
it diagonal, and it is one of the six matrices
$I$, $0$, $T_1$, $T_2$, $B_3=\diag(0,1,1)$,or $B_4=\diag(0,1,0)$. This gives the 20
codifferentials $d_{53}\mcom d_{72}$.

This is not the stable case, and it is not hard to see that the stable case occurs
when $\dim(W)=6$.

\begin{table}[ht]
\begin{center}
\begin{tabular}{ccllrc}
M&W&$\delta$&$\mu$&N&Range\\\hline \\
1&$3$&$\psa{44}4+\psa{33}3+\psa{22}2$&$\psa{11}1$&1&$d_2$\\
1&$3$&$\psa{44}4+\psa{33}3+\psa{22}2$&$0$&5&$d_3\mcom d_7$\\
2&$2$&$\psa{44}4+\psa{33}3$&$\psa{22}1$&2&$d_8,d_9$\\
2&$2$&$\psa{44}4+\psa{33}3$&$0$&25&$d_{10}\mcom d_{34}$\\
3&$1$&$\psa{44}4$&$\psa{31}2+\psa{13}2+\psa{33}1$&2&$d_{35},d_{36}$\\
3&$1$&$\psa{44}4$&$\psa{31}2q+\psa{13}2p+\psa{33}2$&2&$d_{37}(p:q),d_{38}(p:q)$\\
3&$1$&$\psa{44}4$&$\psa{13}2+\psa{33}2$&6&$d_{39}\mcom d_{44}$\\
3&$1$&$\psa{44}4$&$\psa{33}2$&6&$d_{45}\mcom d_{50}$\\
3&$1$&$\psa{44}4$&$\psa{13}2-\psa{31}2$&2&$d_{51}\mcom d_{52}$\\
3&$1$&$\psa{44}4$&$0$&20&$d_{53}\mcom d_{72}$\\
\\ \hline
\end{tabular}
\end{center}
\label{coho04 table}
\caption{Table of Extensions of $\delta$ on $W$ by $\mu$ on $M$}
\end{table}
Note that the simple algebra $d_1$ does not appear in the table above, because it does not arise as an extension.
\section{Hochschild Cohomology and Deformations}
Suppose that $V$ is a vector space, defined over a field $\k$
whose characteristic is not 2 or 3, equipped with an associative multiplication
structure $m:V\tns V\ra V$. The associativity relation can be given in the form
\begin{equation*}
m\circ(m\tns 1)=m\circ(1\tns m).
\end{equation*}

The notion of isomorphism or \emph{equivalence} of associative algebra structures is given as
follows. If $g$ is a linear automorphism of $V$, then define
\begin{equation*}
g^*(m)=g\inv\circ m\circ (g\tns g).
\end{equation*}
Two algebra structures $m$ and $m'$ are equivalent if there is an automorphism
$g$ such that $m'=g^*(m)$. The set of equivalence classes of algebra structures
on $V$ is called the \emph{moduli space} of associative algebras on $V$.

\emph{Hochschild cohomology} was introduced in \cite{hoch}, and was used by Gerstenhaber in \cite{gers} to
classify infinitesimal deformations of associative algebras.

We define the Hochschild coboundary operator $D$ on $\hom(\TV,V)$ by
\begin{align*}
D(\ph)(a_0\mcom a_n)=&a_0\ph(a_1\mcom a_n)+\s{n+1}\ph(a_0
\mcom a_{n-1})a_n\\&+\sum_{i=0}^{n-1}\s{i+1}\ph(a_0\mcom a_{i-1},a_ia_{i+1},a_{i+2}\mcom a_n)
.
\end{align*}
We wish to transform this classical viewpoint into the more modern viewpoint of
associative algebras as being given by codifferentials on a certain coalgebra.
To do this, we first introduce the \emph{parity reversion} $\Pi V$ of a
\zt-graded vector space $V$. If $V=V_e\oplus V_o$ is the decomposition of $V$
into its even and odd parts, then $W=\Pi V$ is the \zt-graded vector space
given by $W_e=V_o$ and $W_o=V_e$. In other words, $W$ is just the space $V$
with the parity of elements reversed.

Given an ordinary associative algebra,
we can view the underlying space $V$ as being \zt-graded, with $V=V_e$.
Then its parity reversion $W$ is again the same space, but now all elements
are considered to be odd. One can avoid this gyration for ordinary spaces,
by introducing a grading by exterior degree on the tensor coalgebra of $V$,
but the idea of parity reversion works equally well when the algebra is \zt-graded,
whereas the method of grading by exterior degree does not.

Denote the tensor (co)-algebra of $W$ by $\TW=\bigoplus_{k=0}^\infty W^k$,
where $W^k$ is the $k$-th tensor power of $W$ and $W^0=\k$. For brevity, the
element in $W^k$ given by the tensor product of the elements $w_i$ in $W$ will
be denoted by $w_1\cdots w_k$. The coalgebra structure on $\TW$ is  given by
\begin{equation*}
\Delta(w_1\cdots w_n)=\sum_{i=0}^n w_1\cdots w_i\tns w_{i+1}\cdots w_n.
\end{equation*}
Define $d:W^2\ra W$ by $d=\pi\circ m\circ (\pi\inv\tns\pi\inv)$, where
$\pi:V\ra W$ is the identity map, which is odd, because it reverses the parity
of elements. Note that $d$ is an odd map. The space $C(W)=\hom(\TW,W)$ is
naturally identifiable with the space of coderivations of $\TW$.  In fact, if
$\ph\in C^k(W)=\hom(W^k,W)$, then $\ph$ is extended to a coderivation of $\TW$
by
\begin{equation*}
\ph(w_1\cdots w_n)=
\sum_{i=0}^{n-k}\s{(w_1\mplus w_i)\ph}w_1\cdots
 w_i\ph(w_{i+1}\cdots w_{i+k})w_{i+k+1}\cdots w_n.
\end{equation*}

The space of coderivations of $\TW$ is equipped with a \zt-graded Lie algebra
structure given by
\begin{equation*}
[\ph,\psi]=\ph\circ\psi-\s{\ph\psi}\psi\circ\ph.
\end{equation*}
The reason that it is more convenient to work with the structure $d$ on $W$
rather than $m$ on $V$ is that the condition of associativity for $m$
translates into the codifferential property $[d,d]=0$.  Moreover, the
Hochschild coboundary operation translates into the coboundary operator $D$ on
$C(W)$, given by
\begin{equation*}
D(\ph)=[d,\ph].
\end{equation*}
This point of view on Hochschild cohomology first appeared in \cite{sta4}.  The
fact that the space of Hochschild cochains is equipped with a graded Lie
algebra structure was noticed much earlier \cite{gers,gers1,gers2,gers3,gers4}.

For notational purposes, we introduce a basis of $C^n(W)$ as follows.  Suppose
that $W=\langle w_1\mcom w_m\rangle$. Then if $I=(i_1\mcom i_n)$ is a
\emph{multi-index}, where $1\le i_k\le m$, denote $w_I=w_{i_1}\cdots w_{i_n}$.
Define $\ph^{I}_i\in C^n(W)$ by
\begin{equation*}
\ph^I_i(w_J)=\delta^I_Jw_i,
\end{equation*}
where $\delta^I_J$ is the Kronecker delta symbol. In order to emphasize the
parity of the element, we will denote $\ph^I_i$ by $\psi^I_i$ when it is an odd
coderivation.

For a multi-index $I=(i_1\mcom i_k)$, denote its \emph{length}  by $\ell(I)=k$. Then
since $W$ is a completely odd space,
the parity of $\ph^I_i$ is given by $|\ph^I_i|=k+1\pmod2$. If
$K$ and $L$ are multi-indices, then denote $KL=(k_1\mcom k_{\ell(K)},l_l\mcom
l_{\ell(L)})$.  Then
\begin{align*}
(\ph^I_i\circ\ph^J_j)(w_K)&=
\sum_{K_1K_2K_3=K}\s{w_{K_1}\ph^J_j} \ph^I_i(w_{K_1},\ph^J_j(w_{K_2}), w_{K_3})
\\&=
\sum_{K_1K_2K_3=K}\s{w_{K_1}\ph^J_j}\delta^I_{K_1jK_3}\delta^J_{K_2}w_i,
\end{align*}
from which it follows that
\begin{equation}\label{braform}
\ph^I_i\circ\ph^J_j=\sum_{k=1}^{\ell(I)}\s{(w_{i_1}\mplus w_{i_{k-1}})\ph^J_j}
\delta^k_j
\ph^{(I,J,k)}_i,
\end{equation}
where $(I,J,k)$ is given by inserting $J$ into $I$ in place of the $k$-th
element of $I$; \ie, $(I,J,k)=(i_1\mcom i_{k-1},j_1\mcom j_{\ell(J)},i_{k+1}\mcom
i_{\ell(I)})$.

Let us explain the notion of an infinitesimal deformation in terms of the
language of coderivations.  We say that
\begin{equation*}
d_t=d+t\psi
\end{equation*}
is an infinitesimal deformation of the codifferential $d$
precisely when $[d_t,d_t]=0 \mod t^2$.
This condition immediately reduces to the cocycle condition $D(\psi)=0$.  Note
that we require $d_t$ to be odd, so that $\psi$ must be an odd coderivation.
One can introduce a more general idea of parameters, allowing both even and odd
parameters, in which case even coderivations play an equal role, but we will
not adopt that point of view in this paper.

For associative algebras, we require that $d$ and $\psi$ lie in $C^2(W)$.
Since in this paper, our algebras are ordinary algebras, so that the
parity of an element in $C^n(W)$ is $n+1$,
elements of $C^2(W)$ are automatically odd.

We need the notion of a versal deformation, in order to
understand how the moduli space is glued together.
To explain versal deformations we introduce
the notion of  a deformation with a local base. For details see
\cite{fi,fi2}.
A local base $A$ is a \zt-graded commutative, unital
$\k$-algebra with an augmentation $\epsilon:A\ra\k$,
whose kernel $\m$ is the unique maximal ideal in $A$,
so that $A$ is a local ring. It follows that $A$ has a
unique decomposition $A=\k\oplus\m$ and $\epsilon$ is
just the projection onto the first factor. Let $W_A=W\tns A$
equipped with the usual structure of a right $A$-module.
Let $T_A(W_A)$ be the tensor algebra of $W_A$ over $A$,
that is $T_A(W_A)=\bigoplus_{k=0}^\infty T^k_A(W_A)$ where $T^0_A(W_A)=A$ and
$T^{k+1}_A(W_A)=T^k(W_A)_A\tns_A W_A$. It is not difficult to show
that $T^k_A(W_A)=T^k(W)\tns A$ in a natural manner, and thus $T_A(W_A)=T(W)\tns A$.

Any $A$-linear map $f:T_A(W)\ra T_A(W)$ is induced by its
restriction to $T(W)\tns \k=T(W)$ so we can view an $A$-linear coderivation
$\delta_A$ on $T_A(W_A)$ as a map $\delta_A:T(W)\ra T(W)\tns A$.
A morphism $f:A\ra B$ induces
a map  $$f_*:\coder_A(T_A(W_A))\ra \coder_B(T_B(W_B))$$ given by
$f_*(\delta_A)=(1\tns f)\delta_A$, moreover if $\delta_A$ is a codifferential
then so is $f_*(A)$.
A codifferential $d_A$ on $T_A(W_A)$ is said to be a deformation
of the codifferential $d$ on $T(W)$ if $\epsilon_*(d_A)=d$.

If $d_A$ is a deformation of $d$ with base $A$ then we can express
\begin{equation*}
 d_A=d+\ph
\end{equation*}
where $\ph:T(W)\ra T(W)\tns\m$. The condition for $d_A$ to be a codifferential is
the Maurer-Cartan equation,
\begin{equation*}
D(\ph)+\frac12[\ph,\ph]=0
\end{equation*}
If $\m^2=0$ we say that $A$ is an infinitesimal algebra and a deformation
with base $A$ is called infinitesimal.

A typical example of an infinitesimal base is $\k[t]/(t^2)$;
 moreover, the classical notion of an infinitesimal deformation:
$d_t=d+t\ph$
is precisely an infinitesimal deformation with base $\k[t]/(t^2)$.

A local algebra $A$ is complete if
\begin{equation*}
 A=\invlim_kA/\m^k
\end{equation*}
A complete, local augmented $\k$-algebra is called formal
and a deformation with a formal base is called a formal deformation,
see \cite{fi2}.
An infinitesimal base is automatically formal, so every infinitesimal deformation
is a formal deformation.

An example of a formal base is $A=\k[[t]]$ and a
deformation of $d$ with base $A$ can be expressed in the form
$$d_t=d+t\psi_1+t^2\psi_2+\dots$$
This is the classical notion of a formal deformation.
It is easy to see that the condition for $d_t$ to be a formal deformation reduces to
\begin{align*}
 D(\psi_{n+1})=-\frac12\sum_{k=1}^{n}[\psi_k,\psi_{n+1-k}],\quad n=0,\dots
\end{align*}

An automorphism of $W_A$ over $A$ is an $A$-linear isomorphism $g_A:W_A\ra W_A$ making
the diagram below commute:

\begin{figure}[h!]
 $$\xymatrix{
 W_A \ar[r]^{g_A} \ar[d]^{\epsilon_*} & W_A \ar[d]^{\epsilon_*} \\
 W \ar[r]^I & W}$$
\end{figure}

The map $g_A$ is induced by its restriction to $T(W)\tns\k$ so we can view $g_A$ as a map
$$g_A:T(W)\ra T(W)\tns A$$
so we ca express $g_A$ in the form
$$g_A=I+\lambda$$
where $\lambda:T(W)\ra T(W)\tns\m$. If $A$ is infinitesimal then $g_A^{-1}=I-\lambda$.

Two deformations $d_A$ and $d_A'$ are said to be equivalent over $A$ if
there is an automorphism $g_A$ of $W_A$ over $A$ such that $g_A^*(d_A)=d_A'$.
In this case we write $d'_A\sim d_A$.

An infinitesimal deformation $d_A$ with base $A$ is called
universal if whenever $d_B$ is an infinitesimal deformation with
base $B$, there is a unique morphism $f:A\ra B$ such that $f_*(d_A)\sim d_B$.

\begin{thm}[\cite{ff2}]
 If $\dim H^2_{odd}(d)<\infty$ then there is a universal
 infinitesimal deformation $\dinf$ of $d$, given by
 $$\dinf=d+\delta^it_i$$
 where $H^2_{odd}(d)=\langle\bar{\delta^i}\rangle$ and $A=\k[t_i]/(t_it_j)$ is
 the base of deformation.
\end{thm}
In the theorem above $\bar\delta^i$ is the cohomology class determined by the cocycle $\delta^i$.

A formal deformation $d_A$ with base $A$ is called versal if given
any formal deformation of $d_B$ with base $B$ there is a morphism
$f:A\ra B$ such that $f_*(d_A)\sim d_B$.
Notice that the difference between the versal and the universal
property of infinitesimal deformations is that $f$ need
not be unique. A versal deformation is called
\emph{miniversal} if $f$ is unique whenever $B$ is
infinitesimal. The basic result about versal deformations is:
\begin{thm}[\cite{fi,fi2,fp1}]
 If $\dim H^2_{odd}(d)<\infty$ then a miniversal deformation of $d$ exists.
\end{thm}

The following result can be used in some special cases to compute the versal deformations.
\begin{thm}
 Suppose $H^2_{odd}(d)=\langle\bar{\delta^i}\rangle$ and
 $[\delta^i,\delta^j]=0$ for all $i,j$ then the infinitesimal deformation
 $$\dinf=d+\delta^it_i$$
 is miniversal, with base $A=\k[[t_i]].$
\end{thm}

The construction of the moduli space as a geometric object
is based on the idea that codifferentials
which can be obtained by deformations with small parameters
are ``close'' to each other. From the small deformations,
we can construct 1-parameter families or even multi-parameter families,
which are defined for small values of the
parameters, except possibly when the parameters vanish.

If $d_t$ is a one parameter family of deformations, then two things can occur.
First, it may happen that
$d_t$ is equivalent to a certain codifferential $d'$ for every small value of $t$ except zero.
Then we say that $d_t$ is a jump deformation from $d$ to $d'$. It will
never occur that $d'$ is equivalent to $d$, so  there are no
jump deformations from a codifferential to itself.
Otherwise, the codifferentials $d_t$ will all be nonequivalent if $t$ is small enough.
In this case, we
say that $d_t$ is a smooth deformation. (In detail, see \cite{fp10}.) 

In \cite{fp10}, it was proved for Lie algebras that given three
codifferentials $d$, $d'$ and $d''$,
if there are jump deformations from $d$ to $d'$ and from $d'$ to $d''$,
then  there is a jump deformation from
$d$ to $d''$. The proof of the corresponding statement for associative algebras
is essentially the same.

Similarly, if there is a jump deformation from $d$ to $d'$, and a family of smooth deformations
$d'_t$, then there is a family $d_t$ of smooth deformations of $d$,
such that every deformation in the image of
$d'_t$ lies in the image of $d_t$, for sufficiently small values of $t$.
In this case, we say that the
smooth deformation of $d$ factors through the jump deformation to $d'$.

In the examples of complex moduli spaces of Lie and associative algebras which we have studied,
it turns out that there is a natural
stratification of the moduli space of $n$-dimensional algebras by orbifolds,
where the codifferentials
on a given strata are connected by smooth deformations
which don't factor through jump deformations.
These smooth deformations determine the local neighborhood structure.

The strata are connected by jump deformations, in the sense that
any smooth deformation from a codifferential on one strata
to another strata factors through a jump deformation.
Moreover, all of the strata are given by projective
orbifolds.

In fact, in all the complex examples we have studied,
the orbifolds either are single points,
or $\C\P^n$ quotiented out by either $\Sigma_{n+1}$ or a subgroup, acting on
$\C\P^n$ by permuting the coordinates.

\section{Deformations of the elements in our moduli space}
We have ordered the codifferentials  so that
a codifferential only deforms to a codifferential earlier on the list.
Partially, this was accomplished the ordering of the different choices of $M$ and $W$.
That such an ordering is possible is due to
the fact that jumps between families have a natural
ordering by descent.

The radical of an algebra $A$ is the same as the radical of its opposite algebra,
ideals in an algebra are the same as the ideals in its opposite algebra $A^\circ$,
and the quotient of the opposite algebra by an ideal is naturally isomorphic to
the quotient of the opposite algebra by the same ideal, it follows that the semisimple quotient
of an algebra is the same as its opposite algebra.
Also the center
of an algebra coincides the center of its opposite algebra.
Moreover, if an algebra $A$ deforms to an
algebra $B$, then its opposite algebra $A^\circ$ deforms to $B^\circ$.
A commutative algebra is isomorphic to its opposite algebra, but an
algebra may be isomorphic to its opposite algebra without being equal to it.
For example, a matrix algebra is always
isomorphic to its opposite algebra, and
the simple $1|1$-dimensional algebra is isomorphic to its opposite,
but neither of these algebras is commutative.

We shall summarize most of the relevant information about the algebras in tables below.
Since there are too many codifferentials to list in a single table,
we will split them up into several tables.  In one set of tables,
we will give the codifferential which represents the algebra, as well as 
information about the cohomology spaces
$H^0$ through $H^3$. In another set of tables, we will note which algebras are pairs of opposite
algebras, give a basis for the center of the algebra, and indicate which algebras it deforms
to. It would take up too much space to give the versal deformations for each of these algebras,
but all of them were computed using the constructive method we have outlined above.
\begin{table}[h]
\begin{center}
\begin{tabular}{lcccc}
Codifferential&$H^0$&$H^1$&$H^2$&$H^3$\\ \hline \\
$d_{1}=\psa{11}1+\psa{12}2+\psa{23}1+\psa{24}2+\psa{32}4+\psa{31}3+\psa{43}3+\psa{44}4$&$1$&$0$&$0$&$0$\\
$d_{2}=\psa{33}3+\psa{44}4+\psa{22}2+\psa{11}1$&$4$&$0$&$0$&$0$\\
$d_{3}=\psa{33}3+\psa{44}4+\psa{22}2+\psa{21}1+\psa{13}1$&$2$&$0$&$0$&$0$\\
$d_{4}=\psa{33}3+\psa{44}4+\psa{22}2+\psa{12}1$&$2$&$0$&$0$&$0$\\
$d_{5}=\psa{33}3+\psa{44}4+\psa{22}2+\psa{21}1$&$2$&$0$&$0$&$0$\\
$d_{6}=\psa{33}3+\psa{44}4+\psa{22}2+\psa{21}1+\psa{12}1$&$4$&$1$&$1$&$1$\\
$d_{7}=\psa{33}3+\psa{44}4+\psa{22}2$&$4$&$1$&$1$&$1$\\
$d_{8}=\psa{33}3+\psa{44}4+\psa{22}1+\psa{31}1+\psa{32}2+\psa{13}1+\psa{23}2$&$4$&$2$&$2$&$2$\\
$d_{9}=\psa{33}3+\psa{44}4+\psa{22}1$&$4$&$2$&$2$&$2$\\
\end{tabular}
\end{center}
\label{d1-d9}
\caption{The cohomology of the algebras $d_1\dots d_9$}
\end{table}

\subsection{The algebras $d_1\dots d_9$}
The algebra $d_1$ represents the matrix algebra $\gl(2,\C)$. As such, it is simple, and so has no ideals,
no deformations, and its center consists of the multiples of the identity, so has dimension 1. Thus
$\dim H^0=1$ and $\dim H^n=0$ otherwise.

The algebra $d_2$ is the semisimple algebra which is the direct sum of four copies of $\C$. Being semisimple,
it is also cohomologically rigid, but it is commutative, so $\dim H^0=4$.

The algebras $d_3$ $d_4$ and $d_5$ are all rigid, with center of dimension 2. The algebras $d_4$ and $d_5$ are opposite
algebras.

The algebra $d_6$ is the direct sum of $C^2$ with the algebra given by adjoining an identity to turn the 1-dimensional trivial algebra into a 2-dimensional unital algebra. It is unital, commutative, is not rigid, and in fact has a jump deformation to $d_2$.

The algebra $d_7$ is the direct sum of the trivial 1-dimensional algebra (which we denote as $\C_0$) with the semisimple
3-dimensional algebra $\C^3$. It is commutative but not unital, and also has a jump deformation to $d_2$.

The algebra $d_8$ arises by as a direct sum of $\C$ and the algebra which arises from adjoining an identity to the 2-dimensional nontrivial nilpotent algebra. It is unital and commutative. It has jump deformations to $d_2$ and $d_6$.

The algebra $d_9$ is the direct sum of the nontrivial 2-dimensional nilpotent algebra and $\C^2$. It is not unital, but is commutative. It has deformations to $d_2$, $d_6$ and $d_7$.

\begin{table}[h,t]
\begin{center}
\begin{tabular}{lcccc}
Codifferential&$H^0$&$H^1$&$H^2$&$H^3$\\ \hline \\
$d_{10}=\psa{33}3+\psa{44}4+\psa{31}1+\psa{32}2+\psa{14}1$&$0$&$0$&$0$&$0$\\
$d_{11}=\psa{33}3+\psa{44}4+\psa{41}1+\psa{13}1+\psa{23}2$&$0$&$0$&$0$&$0$\\
$d_{12}=\psa{33}3+\psa{44}4+\psa{31}1+\psa{42}2$&$0$&$0$&$0$&$0$\\
$d_{13}=\psa{33}3+\psa{44}4+\psa{13}1+\psa{24}2$&$0$&$0$&$0$&$0$\\
$d_{14}=\psa{33}3+\psa{44}4+\psa{31}1+\psa{24}2$&$0$&$0$&$0$&$0$\\
$d_{15}=\psa{33}3+\psa{44}4+\psa{31}1+\psa{14}1+\psa{42}2$&$0$&$0$&$0$&$0$\\
$d_{16}=\psa{33}3+\psa{44}4+\psa{13}1+\psa{41}1+\psa{24}2$&$0$&$0$&$0$&$0$\\
$d_{17}=\psa{33}3+\psa{44}4+\psa{31}1+\psa{32}2+\psa{14}1+\psa{24}2$&$1$&$3$&$0$&$0$\\
$d_{18}=\psa{33}3+\psa{44}4+\psa{31}1+\psa{32}2$&$1$&$3$&$0$&$0$\\
$d_{19}=\psa{33}3+\psa{44}4+\psa{13}1+\psa{23}2$&$1$&$3$&$0$&$0$\\
$d_{20}=\psa{33}3+\psa{44}4+\psa{23}2+\psa{31}1$&$1$&$1$&$0$&$1$\\
$d_{21}=\psa{33}3+\psa{44}4+\psa{13}1+\psa{32}2+\psa{41}1+\psa{24}2$&$1$&$1$&$1$&$1$\\
$d_{22}=\psa{33}3+\psa{44}4+\psa{32}2+\psa{23}2+\psa{31}1$&$2$&$1$&$1$&$1$\\
$d_{23}=\psa{33}3+\psa{44}4+\psa{32}2+\psa{23}2+\psa{13}1$&$2$&$1$&$1$&$1$\\
$d_{24}=\psa{33}3+\psa{44}4+\psa{31}1+\psa{14}1$&$2$&$1$&$1$&$1$\\
$d_{25}=\psa{33}3+\psa{44}4+\psa{31}1+\psa{42}2+\psa{24}2$&$2$&$1$&$1$&$1$\\
$d_{26}=\psa{33}3+\psa{44}4+\psa{13}1+\psa{42}2+\psa{24}2$&$2$&$1$&$1$&$1$\\
$d_{27}=\psa{33}3+\psa{44}4+\psa{13}1+\psa{41}1+\psa{42}2+\psa{24}2$&$2$&$1$&$1$&$1$\\
$d_{28}=\psa{33}3+\psa{44}4+\psa{31}1+\psa{14}1+\psa{42}2+\psa{24}2$&$2$&$1$&$1$&$1$\\
$d_{29}=\psa{33}3+\psa{44}4+\psa{32}2$&$2$&$1$&$2$&$2$\\
$d_{30}=\psa{33}3+\psa{44}4+\psa{23}2$&$2$&$1$&$2$&$2$\\
$d_{31}=\psa{33}3+\psa{44}4+\psa{32}2+\psa{23}2$&$4$&$2$&$2$&$2$\\
$d_{32}=\psa{33}3+\psa{44}4+\psa{31}1+\psa{13}1+\psa{42}2+\psa{24}2$&$4$&$2$&$2$&$2$\\
$d_{33}=\psa{33}3+\psa{44}4+\psa{31}1+\psa{13}1+\psa{32}2+\psa{23}2$&$4$&$4$&$6$&$12$\\
$d_{34}=\psa{33}3+\psa{44}4$&$4$&$4$&$8$&$16$\\
\end{tabular}
\end{center}
\label{d10-d34}
\caption{The cohomology of the algebras $d_{10}\dots d_{34}$}
\end{table}

\subsection{The algebras $d_{10}\dots d_{34}$}
The algebras $d_{10}\dots d_{20}$ are all nonunital, noncommutative, and rigid. The pairs of opposite algebras
are $d_{10}$ and $d_{11}$, $d_{12}$ and $d_{13}$, $d_{15}$ and $d_{16}$, and $d_{18}$ and $d_{19}$. The algebras
$d_{14}$, $d_{17}$ and $d_{20}$ are all isomorphic to their opposite algebras.

The algebra $d_{21}$ is unital, but not commutative, and it has a jump deformation to $d_1$.

The algebras $d_{22}$ and $d_{23}$ are nonunital, noncommutative opposite algebras, with $d_{22}$ having a jump
deformation to $d_{5}$, while $d_{23}$ jumps to its opposite algebra $d_{4}$. Similarly $d_{25}$ and $d_{26}$ are also
nonunital, noncommutative opposite algebras which jump to the same two elements in the same order.

The algebra
$d_{24}$ is a nonunital, noncommutative algebra which is isomorphic to its opposite algebra, and it jumps to $d_3$.
The algebras $d_{27}$ and $d_{28}$ are unital, noncommutative opposite algebras both of which have jump deformations
to $d_3$.

The algebras $d_{29}$ and $d_{30}$ are nonunital, nonunital opposite algebras, with $d_{29}$ jumping  to $d_3$ and $d_5$,
while $d_{30}$ jumps to $d_3$ and $d_4$.

The algebra $d_{31}$ is nonunital but is commutative, and it has jump deformations to $d_2$, $d_6$ and $d_7$, all of which
are commutative. Note that a commutative algebra may deform to a noncommutative algebra, but the converse is impossible.

The algebra $d_{32}$ is both unital and commutative and jumps to $d_2$ and $d_6$. Note that a unital algebra can only deform
to another unital algebra, and both $d_2$ and $d_6$ are unital.

The algebra $d_{33}$ which arises by first taking the trivial 2-dimensional algebra, adding a multiplicative identity to make
it unital, and then taking a direct sum with $\C$, is both unital and commutative. It deforms to $d_2$, $d_3$, $d_6$ and $d_8$.
Note that $d_3$ is not commutative, illustrating the fact that a commutative algebra can deform to a noncommutative algebra.

The final algebra in this group, $d_{34}$, is the direct sum of the trivial 2-dimensional algebra $\C_0^2$ with $\C^2$, so it
is not unital, but is commutative. This algebra has a lot of deformations, with jump deformations to $d_2$, $d_4$, $d_5$, $d_6$,
$d_7$ and $d_9$. Note that even though $d_{34}$ is isomorphic to its opposite, it has jump deformations to $d_4$ and $d_5$, which
are not their own opposites. However, they are opposite algebras, illustrating the fact that if an algebra which is isomorphic
to its opposite deforms to another algebra, it also deforms to the opposite of that algebra.
\begin{table}[h,t]
\begin{center}
\begin{tabular}{lcccc}
Codifferential&$H^0$&$H^1$&$H^2$&$H^3$\\ \hline \\
$d_{35}=\psa{44}4+\psa{31}2+\psa{13}2+\psa{33}1+\psa{41}1+\psa{42}2+\psa{43}3+\psa{14}1+\psa{24}2+\psa{34}3$&$4$&$3$&$3$&$3$\\
$d_{36}=\psa{44}4+\psa{31}2+\psa{13}2+\psa{33}1$&$4$&$3$&$3$&$3$\\
$d_{37}(p:q)=\psa{44}4+q\psa{31}2+p\psa{13}2+\psa{33}2+\psa{41}1+\psa{42}2+\psa{43}3+\psa{14}1+\psa{24}2+\psa{34}3$&$2$&$2$&$1$&$0$\\
$d_{37}(1:1)=\psa{44}4+\psa{31}2+\psa{13}2+\psa{33}2+\psa{41}1+\psa{42}2+\psa{43}3+\psa{14}1+\psa{24}2+\psa{34}3$&$4$&$4$&$5$&$6$\\
$d_{37}(1:-1)=\psa{44}4-\psa{31}2+\psa{13}2+\psa{33}2+\psa{41}1+\psa{42}2+\psa{43}3+\psa{14}1+\psa{24}2+\psa{34}3$&$2$&$2$&$1$&$1$\\
$d_{37}(1:0)=\psa{44}4+\psa{13}2+\psa{33}2+\psa{41}1+\psa{42}2+\psa{43}3+\psa{14}1+\psa{24}2+\psa{34}3$&$2$&$2$&$3$&$5$\\
$d_{37}(0:1)=\psa{44}4+\psa{31}2+\psa{33}2+\psa{41}1+\psa{42}2+\psa{43}3+\psa{14}1+\psa{24}2+\psa{34}3$&$2$&$2$&$3$&$5$\\
$d_{37}(0:0)=\psa{44}4+\psa{33}2+\psa{41}1+\psa{42}2+\psa{43}3+\psa{14}1+\psa{24}2+\psa{34}3$&$4$&$5$&$7$&$13$\\
$d_{38}(p:q)=\psa{44}4+q\psa{31}2+p\psa{13}2+\psa{33}2$&$2$&$2$&$3$&$3$\\
$d_{38}(1:1)=\psa{44}4+\psa{31}2+\psa{13}2+\psa{33}2$&$4$&$4$&$5$&$7$\\
$d_{38}(1:-1)=\psa{44}4-\psa{31}2+\psa{13}2+\psa{33}2$&$2$&$2$&$3$&$4$\\
$d_{38}(1:0)=\psa{44}4+\psa{13}2+\psa{33}2$&$2$&$2$&$5$&$8$\\
$d_{38}(0:1)=\psa{44}4+\psa{31}2+\psa{33}2$&$2$&$2$&$5$&$8$\\
$d_{38}(0:0)=\psa{44}4+\psa{33}2$&$4$&$5$&$9$&$17$\\
$d_{39}=\psa{44}4+\psa{13}2+\psa{33}2+\psa{41}1+\psa{42}2+\psa{43}1-\psa{34}1+\psa{24}2+\psa{34}3$&$1$&$1$&$0$&$0$\\
$d_{40}=\psa{44}4+\psa{13}2+\psa{33}2-\psa{43}1+\psa{43}3+\psa{14}1+\psa{34}1$&$1$&$1$&$1$&$0$\\
$d_{41}=\psa{44}4+\psa{13}2+\psa{33}2+\psa{41}1+\psa{42}2+\psa{43}1$&$0$&$1$&$1$&$1$\\
$d_{42}=\psa{44}4+\psa{13}2+\psa{33}2-\psa{34}1+\psa{24}2+\psa{34}3$&$0$&$1$&$1$&$1$\\
$d_{43}=\psa{44}4+\psa{13}2+\psa{33}2-\psa{43}1+\psa{43}3+\psa{14}1+\psa{24}2+\psa{34}3$&$0$&$1$&$1$&$1$\\
$d_{44}=\psa{44}4+\psa{13}2+\psa{33}2+\psa{41}1+\psa{42}2+\psa{43}3+\psa{14}1+\psa{34}1$&$0$&$1$&$1$&$1$\\
$d_{45}=\psa{44}4+\psa{33}2+\psa{41}1+\psa{24}2+\psa{42}2+\psa{34}3+\psa{43}3$&$2$&$2$&$2$&$2$\\
$d_{46}=\psa{44}4+\psa{33}2+\psa{14}1+\psa{24}2+\psa{42}2+\psa{34}3+\psa{43}3$&$2$&$2$&$2$&$2$\\
$d_{47}=\psa{44}4+\psa{33}2+\psa{14}1$&$2$&$2$&$3$&$3$\\
$d_{48}=\psa{44}4+\psa{33}2+\psa{41}1$&$2$&$2$&$3$&$3$\\
$d_{49}=\psa{44}4+\psa{33}2+\psa{14}1+\psa{41}1$&$4$&$3$&$3$&$3$\\
$d_{50}=\psa{44}4+\psa{33}2+\psa{24}2+\psa{42}2+\psa{34}3+\psa{43}3$&$4$&$3$&$3$&$3$\\
$d_{51}=\psa{44}4-\psa{31}2+\psa{13}2+\psa{41}1+\psa{42}2+\psa{43}3+\psa{14}1+\psa{24}2+\psa{34}3$&$2$&$4$&$6$&$8$\\
$d_{52}=\psa{44}4-\psa{31}2+\psa{13}2$&$2$&$4$&$5$&$8$\\
\end{tabular}
\end{center}
\label{d35-d52}
\caption{The cohomology of the algebras $d_{35}\dots d_{52}$}
\end{table}

\subsection{The algebras $d_{35}\dots d_{52}$}
The algebra $d_{35}$ which arises by adjoining an identity to  the 3-dimensional nilpotent algebra
$d_{19}=\psa{31}2+\psa{13}2+\psa{33}1$ is both unital and commutative. It has jump deformations to
$d_{2}$, $d_6$, $d_8$ and $d_{32}$.

The algebra $d_{36}$ which is the direct sum of $\C$ and the 3-dimensional nilpotent algebra $d_{19}$ above, is nonunital
but commutative, and it deforms to $d_2$ $d_6$, $d_7$, $d_8$, $d_9$, and $d_{31}$.

The family of algebras $d_{37}(p:q)$ is parameterized by $\P^1/\Sigma_2$, which means that it is a projective family,
in the sense that $d_{37}(p:q)\sim d_{37}(up:uq)$ when $u\in\C^*$, and is invariant under the action of $\Sigma_2$ by
interchanging coordinates, in the sense that $d_{37}(p:q)\sim d_{37}(q:p)$.

We remark that there is an element
$d_{37}(0:0)$ corresponding to what is called the generic point in $\P^1$. This point is usually omitted in the definition of $\P^1$, because including this generic point makes $\P^1$ a non-Hausdorff space. In fact, this non-Hausdorff behavior is reflected in
the deformations of the point $d_{37}(0:0)$, so the inclusion of the corresponding codifferential in the family here is quite natural.
With the families, there is a generic deformation pattern, and then there are some special values of the parameter $(p:q)$ for which the deformation pattern is not generic in the sense that there are additional deformations.  Generically, this family consists
of unital but not commutative algebras. Because the opposite algebra to $d_{37}(p:q)$ is $d_{37}(q:p)$ which is isomorphic to the
original algebra, all of the elements of this family are isomorphic to their own opposite algebras.

Generically, an element in
this algebra deforms in a smooth way to other elements in the family, and these are the only deformations. We say that the
deformations are along the family. In fact, in every family of codifferentials, there are always smooth deformations along the family. In this case, these are the only deformations which occur generically.

The element $d_{37}(1:0)$ has additional jump deformations to $d_3$, $d_{27}$ and $d_{28}$. The element $d_{1:1}$ is commutative,
and has
additional jump deformations to $d_2$, $d_6$ $d_8$, $d_{32}$ and $d_{35}$, all of which are unital, commutative algebras.

If an element in a family has a deformation to an algebra, then the generic element in the family will also deform to it.
Moreover, the generic element always has jump deformations to all other
 elements in the family, so $d_{37}(0:0)$ has
jump deformations to $d_{37}(p:q)$ for all $(p:q)$ except $(0:0)$.
Thus we automatically know that $d_{37}(0:0)$ has jump deformations the elements to which $d_{37}(1:0)$ and $d_{37}(1:1)$ deform.
In addition, there is a jump deformation from $d_{37}(0:0)$ to $d_{33}$. We also note that $d_{37}(0:0)$ is commutative.

The family $d_{38}(p:q)$ is also parameterized projectively by $\P^1/\Sigma_2$. Generically, the elements of the family are not commutative, and the only deformations are smooth deformations along the family.

The algebra $d_{38}(1:0)$ also has jump deformations to $d_3$, $d_4$, $d_5$, $d_{22}$, $d_{23}$, $d_{29}$, and $d_{30}$.
The algebra $d_{38}(1:1)$ is commutative and has additional jump deformations to $d_2$, $d_6$, $d_7$, $d_8$, $d_9$, $d_{31}$, and $d_{36}$. Finally, the generic element, which is commutative, in addition to all deformations above, and jump deformations to every
other element of the family, also has jumps to $d_3$, $d_{33}$ and $d_{34}$.

The algebras $d_{39}\dots d_{48}$ are neither unital nor commutative. The algebras $d_{39}$ are isomorphic to their opposite algebras. The algebra $d_{39}$ is rigid, while the algebra $d_{40}$ has a jump deformation to $d_1$. The algebras $d_{41}$ and $d_{42}$ are opposite algebras, with $d_{41}$ having a jump deformation to $d_{10}$, while $d_{42}$ deforms to the opposite algebra $d_{11}$.
The algebras $d_{43}$ and $d_{44}$ are opposites, with $d_{43}$ jumping to $d_{13}$ and $d_{44}$ jumping to its opposite $d_{12}$,

The algebras $d_{45}$ and $d_{46}$ are opposites, with $d_{45}$ jumping to $d_5$, $d_{22}$ and $d_{25}$, while $d_{46}$ jumps to
their opposite algebras $d_{4}$, $d_{23}$ and $d{26}$. The algebras $d_{47}$ and $d_{48}$ are opposite algebras, with
$d_{47}$ having jump deformations to $d_3$, $d_4$, $d_{24}$, $d_{26}$, $d_{27}$, and  $d_{30}$, while $d_{48}$ jumps to
$d_3$, $d_5$, $d_{24}$, $d_{25}$, $d_{28}$, and $d_{29}$.

The algebras $d_{49}$ and $d_{50}$ are nonunital but are commutative, with $d_{49}$ jumping to $d_2$, $d_6$, $d_7$, $d_9$, $d_{31}$, and $d_{32}$, while $d_{50}$ jumps to $d_2$, $d_6$, $d_7$, $d_8$, and $d_{31}$.

The algebra $d_{51}$ is unital, with jump deformations to $d_1$, $d_{21}$, $d_{37}(1:-1)$, and deforms smoothly near $d_{37}(1:-1)$. This type of smooth deformation is said to factor through the jump deformation to $d_{37}(1:-1)$.

The algebra $d_{52}$ is neither unital nor commutative, but is isomorphic to its own opposite algebra. It has jump deformations to
$d_{20}$ and $d_{38}(1: -1)$, as well as deforming smoothly in a neighborhood of $d_{38}(1: -1)$.

\begin{table}[h,t]
\begin{center}
\begin{tabular}{lcccc}
Codifferential&$H^0$&$H^1$&$H^2$&$H^3$\\ \hline \\
$d_{53}=\psa{44}4+\psa{14}1+\psa{24}2+\psa{43}3$&$0$&$4$&$0$&$8$\\
$d_{54}=\psa{44}4+\psa{41}1+\psa{42}2+\psa{34}3$&$0$&$4$&$0$&$8$\\
$d_{55}=\psa{44}4+\psa{41}1+\psa{42}2+\psa{43}3$&$0$&$8$&$0$&$0$\\
$d_{56}=\psa{44}4+\psa{14}1+\psa{24}2+\psa{34}3$&$0$&$8$&$0$&$0$\\
$d_{57}=\psa{44}4+\psa{14}1+\psa{41}1+\psa{24}2+\psa{43}3$&$1$&$2$&$3$&$5$\\
$d_{58}=\psa{44}4+\psa{14}1+\psa{41}1+\psa{42}2$&$2$&$2$&$3$&$4$\\
$d_{59}=\psa{44}4+\psa{14}1+\psa{41}1+\psa{24}2$&$2$&$2$&$3$&$4$\\
$d_{60}=\psa{44}4+\psa{41}1+\psa{34}3$&$1$&$2$&$4$&$6$\\
$d_{61}=\psa{44}4+\psa{14}1+\psa{41}1+\psa{42}2+\psa{43}3$&$1$&$4$&$4$&$4$\\
$d_{62}=\psa{44}4+\psa{14}1+\psa{41}1+\psa{24}2+\psa{34}3$&$1$&$4$&$4$&$4$\\
$d_{63}=\psa{44}4+\psa{41}1+\psa{42}2$&$1$&$4$&$5$&$5$\\
$d_{64}=\psa{44}4+\psa{14}1+\psa{24}2$&$1$&$4$&$5$&$5$\\
$d_{65}=\psa{44}4+\psa{14}1+\psa{41}1+\psa{24}2+\psa{42}2$&$4$&$5$&$7$&$13$\\
$d_{66}=\psa{44}4+\psa{14}1+\psa{41}1+\psa{24}2+\psa{42}2+\psa{43}3$&$2$&$4$&$8$&$16$\\
$d_{67}=\psa{44}4+\psa{14}1+\psa{41}1+\psa{24}2+\psa{42}2+\psa{34}3$&$2$&$4$&$8$&$16$\\
$d_{68}=\psa{44}4+\psa{14}1+\psa{41}1$&$4$&$5$&$9$&$17$\\
$d_{69}=\psa{44}4+\psa{41}1$&$2$&$4$&$10$&$20$\\
$d_{70}=\psa{44}4+\psa{14}1$&$2$&$4$&$10$&$20$\\
$d_{71}=\psa{44}4+\psa{14}1+\psa{41}1+\psa{24}2+\psa{42}2+\psa{34}3+\psa{43}3$&$4$&$9$&$24$&$72$\\
$d_{72}=\psa{44}4$&$4$&$9$&$27$&$81$\\
\end{tabular}
\end{center}
\label{d53-d72}
\caption{The cohomology of the algebras $d_{53}\dots d_{72}$}
\end{table}

\subsection{The algebras $d_{53}\dots d_{72}$}
The algebras $d_{53}\dots d_{68}$ are neither unital nor commutative. The algebras $d_{53}$ and $d_{54}$, and $d_{55}$ and $d_{56}$,
are pairs of opposite algebras, and they are all rigid.
The algebra $d_{57}$ is isomorphic to its opposite algebra, and jumps to $d_{14}$, $d_{20}$ and  $d_{39}$.
the algebras $d_{58}$ and $d_{59}$ are opposites, with the former jumping to $d_3$, $d_5$, $d_{22}$, $d_{27}$, and $d_{29}$,
while the latter jumps to $d_3$, $d_4$, $d_{23}$, $d_{28}$, and $d_{30}$.

The algebra $d_{60}$ is isomorphic to its opposite algebra, and jumps to $d_1$, $d_{15}$, $d_{16}$, $d_{20}$, $d_{21}$, and $d_{40}$.
The algebras $d_{61}$ and $d_{62}$ are opposites, with $d_{61}$ jumping to $d_{12}$, $d_{18}$ and  $d_{44}$, while
$d_{62}$ jumps to $d_{13}$ $d_{19}$ and $d_{43}$. The algebras $d_{63}$ and $d_{64}$ are opposite algebras, with
$d_{63}$ having jump deformations to $d_{10}$, $d_{17}$, $d_{18}$, and $d_{41}$, while $d_{64}$ jumps to
$d_{11}$, $d_{17}$, $d_{19}$, and $d_{42}$.

The algebra $d_{65}$ is its own opposite algebra, and it jumps to $d_2$, $d_3$, $d_6$, $d_7$, $d_8$, $d_{24}$, $d_{31}$, $d_{33}$,
and $d_{50}$. The algebras $d_{66}$ and $d_{67}$ are opposite algebras, with $d_{66}$ having jump deformations to
$d_5$, $d_{10}$, $d_{15}$, $d_{22}$, $d_{25}$, and $d_{45}$, and $d_{67}$ jumping to
$d_4$, $d_{11}$, $d_{16}$, $d_{23}$, $d_{26}$, and  $d_{46}$.

The algebra $d_{68}$ is not unital, but is commutative, and it jumps to $d_2$, $d_4$, $d_5$, $d_6$, $d_7$,
$d_9$, $d_{25}$, $d_{26}$, $d_{31}$, $d_{32}$, $d_{34}$ and $d_{49}$. The opposite algebras $d_{69}$ and $d_{70}$ are neither unital
nor commutative, with the former deforming to $d_3$, $d_5$, $d_{11}$, $d_{14}$, $d_{15}$, $d_{24}$, $d_{25}$, $d_{28}$,
$d_{29}$, $d_{48}$, and the latter to $d_3$, $d_4$, $d_{10}$, $d_{13}$, $d_{14}$, $d_{16}$, $d_{24}$, $d_{26}$, $d_{27}$, $d_{30}$, and $d_{47}$.

The algebra $d_{71}$ is  the algebra arising by adjoining an
identity to the trivial 3-dimensional algebra $\C^3_0$, so it is unital and commutative, and it has deformations to every unital
algebra, that is, to $d_1$, $d_2$, $d_3$, $d_6$, $d_8$, $d_{17}$, $d_{21}$,
 $d_{27}$, $d_{28}$, $d_{32}$, $d_{33}$, $d_{35}$, every element of the family $d_{37}(p:q)$, and $d_{51}$.

Finally, the algebra $d_{72}$ is the direct sum of $\C$ and $\C^3_0$, so it is nonunital and commutative. It has jump deformations to
$d_2\dots d_9$, $d_{18}$, $d_{19}$, $d_{20}$, $d_{22}$, $d_{23}$, $d_{29}$, $d_{30}$, $d_{31}$, $d_{33}$,
$d_{34}$, $d_{36}$, all members of the family  $d_{38}(p:q)$, and $d_{52}$.

\begin{table}[h,t]
\begin{center}
\begin{tabular}{lll}
Codifferential&Gabriel Number&Structure\\ \hline \\
$d_{1}$&10& $\gl(2,\C)$\\
$d_2$&1&$\C^4$\\
$d_3$&13&$\C\oplus\{\left[\begin{smallmatrix}a&b\\0&c\end{smallmatrix}\right]|a,b,c\in\C\}$\\
$d_6$&2&$\C^2\oplus\C[x]/(x^2)$\\
$d_8$&4&$\C\oplus\C[x]/(x^3)$\\
$d_{17}$&17&$\left\{\left[\begin{smallmatrix}a&0&0\\0&a&0\\c&d&b\end{smallmatrix}\right]|a,b,c,d\in\C\right\}$\\
$d_{21}$&11&$\left\{\left[\begin{smallmatrix}a&0&0&0\\0&a&0&d\\c&0&b&0\\0&0&0&b\end{smallmatrix}\right]|a,b,c,d\in\C\right\}$\\
$d_{27}$&15&$\left\{\left[\begin{smallmatrix}a&c&d\\0&a&0\\0&0&b\end{smallmatrix}\right]|a,b,c,d\in\C\right\}$\\
$d_{28}$&14&$\left\{\left[\begin{smallmatrix}a&0&0\\c&a&0\\d&0&b\end{smallmatrix}\right]|a,b,c,d\in\C\right\}$\\
$d_{32}$&3&$\C[x]/(x^2)\oplus\C[y]/(y^2)$\\
$d_{33}$&6&$\C\oplus\C[x,y]/(x^2,xy,y^2)$\\
$d_{35}$&5&$\C[x]/(x^4)$\\
$d_{37}(p:q)$&$18(u)$&$\C[x,y]/(x^2,y^2,yx,-uxy),\quad u\ne-1$\\
$d_{37}(1:-1)$&19&$\C[x,y]/(y^2,x^2+yx,xy+yx)$\\
$d_{37}(1:1)$&7&$\C[x,y]/(x^2,y^2)$\\
$d_{37}(1:0)$&16&$\C[x,y]/(x^2,y^2,yx)$\\
$d_{37}(0:0)$&8&$\C[x,y]/(x^3,xy,y^2)$\\
$d_{51}$&12&$\bigwedge\C^2$\\
$d_{71}$&9&$\C[x,y,z]/(x,y,z)^2$\\
\end{tabular}
\end{center}
\label{Structure of Unital 4-dimensional algebras}
\caption{The structure of 4-dimensional unital algebras}
\end{table}

\section{Unital Algebras}
There are 15 unital algebras, including all of the elements in the family $d_{37}(p:q)$.
According to \cite{mas}, unital complex 4-dimensional associative algebras were classified by P. Gabriel \cite{gab}, and
this classification is in agreement with the unital algebras we determined by our methods. Note that no nilpotent algebra
can be unital, so the classification of the nonnilpotent algebras given here is sufficient to determine all of the unital
algebras.

\begin{table}[h,t]
\begin{center}
\begin{tabular}{ll}
Codifferential&Structure\\ \hline \\
$d_{2}$& $\C^4$\\
$d_{6}$& $\C^2\oplus\C[x]/(x^2)$\\
$d_7$&$\C^3\oplus\C_0$\\
$d_8$&$\C\oplus\C[x]/(x^3)$\\
$d_9$&$\C^2\oplus x\C[x]/(x^3)$\\
$d_{31}$&$\C\oplus\C_0\oplus\C[x]/(x^2)$\\
$d_{32}$&$\C[x]/(x^2)\oplus\C[y]/(y^2)$\\
$d_{33}$&$\C\oplus\C[x,y]/(x^2,xy,y^2)$\\
$d_{34}$&$\C^2\oplus\C^2_0$\\
$d_{35}$&$\C[x]/(x^4)$\\
$d_{36}$&$\C\oplus x\C[x]/(x^4)$\\
$d_{37}(1:1)$&$\C[x,y]/(x^2,y^2)$\\
$d_{37}(0:0)$&$\C[x,y]/(x^3,xy,y^2)$\\
$d_{38}(1:1)$&$\C\oplus(x,y)\le\C\oplus\C[x,y]/(x^2-xy,y^2)$\\
$d_{38}(0:0)$&$\C\oplus\C_0\oplus x\C[x]/(x^3)$\\
$d_{49}$&$\C[x]/(x^2)\oplus y\C[y]/(y^3)$\\
$d_{50}$&$\C_0\oplus\C[x]/(x^3)$\\
$d_{68}$&$\C_0^2\oplus C[x]/(x^2)$\\
$d_{71}$&$\C[x,y,z]/(x,y,z)^2$\\
$d_{72}$&$\C\oplus\C_0^3$\\
\end{tabular}
\end{center}
\label{Commutative Nonnilpotent}
\caption{The structure of nonnilpotent 4-dimensional commutative algebras}
\end{table}
\bigskip

\begin{figure}[htp]%
\begin{center}
\includegraphics[scale=0.42]{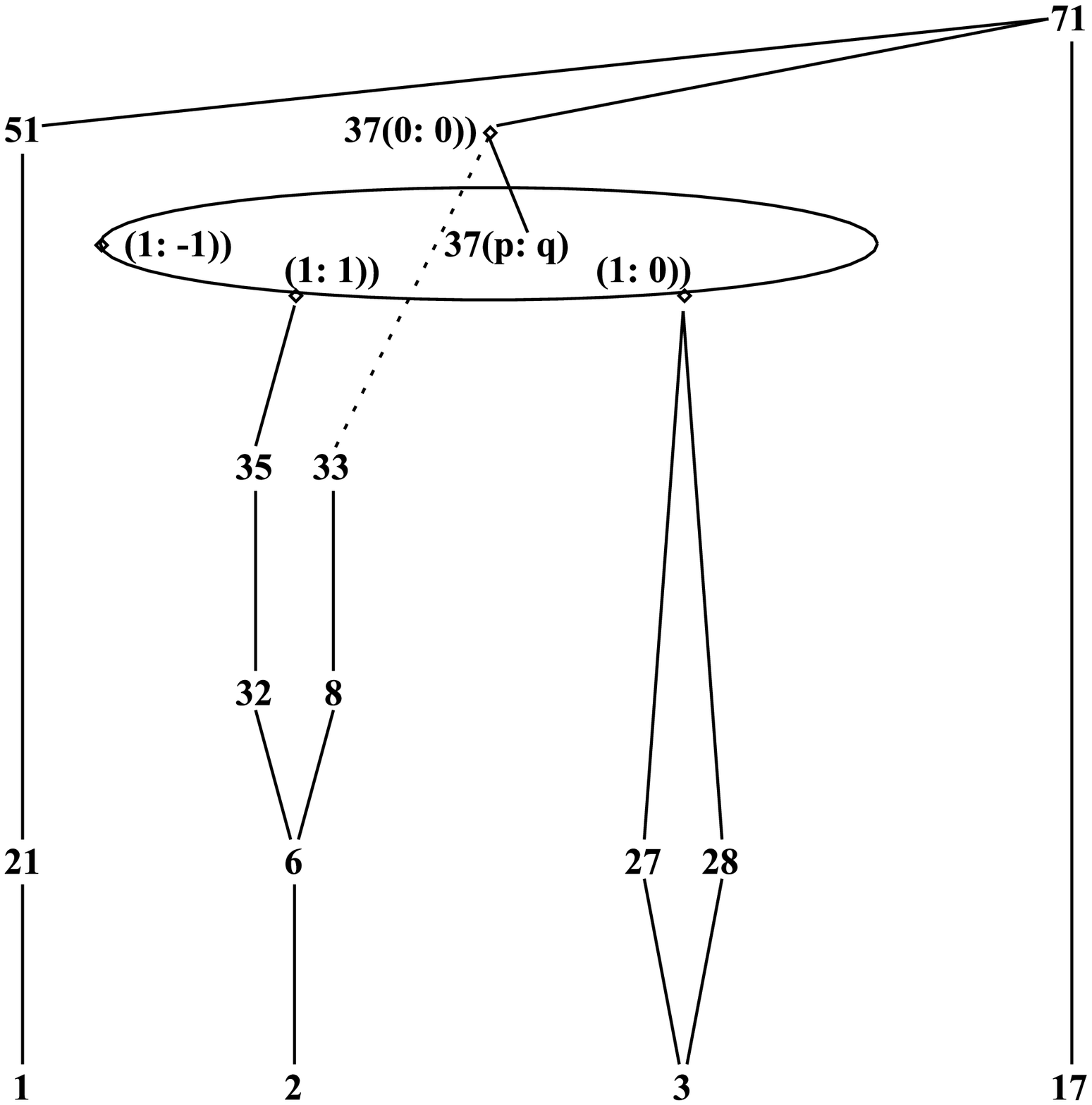}
\caption{Deformations between unital algebras}%
\label{unital}%
\end{center}
\end{figure}

\section{Commutative Algebras}
There are 20 distinct nonnilpotent commutative algebras, of which 9 are unital. Every commutative algebra is a direct sum of
algebras which are ideals in quotients of polynomial algebras. Every finite dimensional unital commutative algebra is a quotient of a polynomial algebra, while every
finite dimensional nonunital algebra is an ideal in such an algebra.  The
algebra $\C$ is representable as $\C[x]/(x)$, while the trivial algebra $\C_0$ is representable as the ideal $x\C[x]/(x^2)$.
In Table 7, the ideal $(x,y)$ in $\C[x,y]/(x^2-xy,y^2)$ has dimension 3 as a vector space over $\C$, and the algebra
$d_{38}(1:1)$ is expressed as a direct sum of $\C$ and that ideal, which gives a 4-dimensional algebra.

For completeness here, in Table 8, we give the nilpotent commutative algebras as well. The codifferential number
given relates to the dsscription of codifferentials which will appear in a sequel. These algebras were
classified by Hazlett \cite{hazl}, and also given in \cite{mazz}. There are 8 nontrivial commutative algebras.

We note that commutative algebras may deform into noncommutative algebras, but noncommutative algebras never deform into
a commutative algebras.  The fact that commutative algebras have noncommutative deformations plays an important role in physics,
and deformation quantization describes a certain type of deformation of a commutative algebra into a noncommutative one.

\begin{table}[h,t]
\begin{center}
\begin{tabular}{ll}
Codifferential&Structure\\ \hline \\
$d_{74}$&$x\C[x]/(x^5)$\\
$d_{75}(1:1)$&$(x,y)\le \C[x,y]/(x^2-y^2,yx^2,xy^2) $\\
$d_{75}(0:0)=d_{86}(1:1)$&$\C_0\oplus(x,y)\le\C_0\oplus\C[x,y]/(x^2-y^2,xy)$\\
$d_{76}$&$(x,y)\le\C[x,y]/(y^3-x^2,xy)$\\
$d_{79}(1:1)$&$(x,y)\le\C[x,y]/(y^2,x^2y,x^3)$\\
$d_{79}(0:0)=d_{86}(0:0)$&$\C_0^2\oplus x\C[x]/(x^3)$\\
$d_{83}$&$\C_0\oplus x\C[x]/(x^4)$\\
$d_{85}$&$(x,y,z)\le\C[x,y,z]/(x^2-y^2,y^2-yz,xy,xz,z^2)$\\
$d_0$&$\C_0^4$\\
\end{tabular}
\end{center}
\label{Commutative Nilpotent}
\caption{Nilpotent 4-dimensional commutative algebras}
\end{table}
\bigskip

\begin{figure}[htp]%
\begin{center}
\includegraphics[scale=0.45]{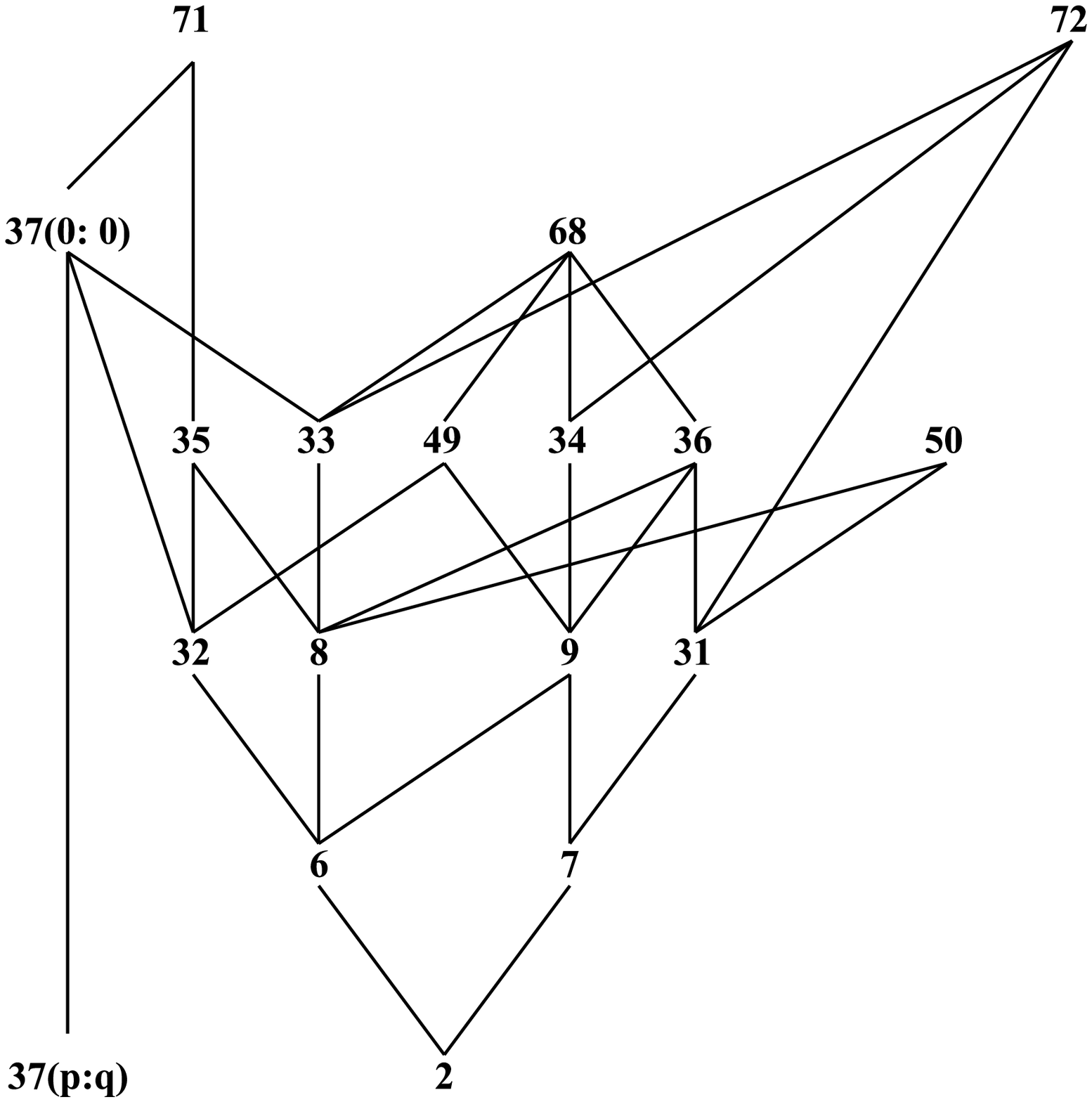}
\caption{Deformations between nonnilpotent commutative algebras}%
\label{commutative}%
\end{center}
\end{figure}

\section{Levels of algebras}
It would be difficult to construct a picture showing the jump
deformations for all 72 families of nonnilpotent complex 4-dimensional
algebras, as we did for the unital and commutative algebras, because
there are too many of them.  Instead, we give a table showing the
levels of each algebra. To define the level, we say that a rigid
algebra has level 1, an algebra which has only jump deformations to an
algebra on level one has level two and so on. To be on level $k+1$, an
algebra must have a jump deformation to an algebra on level $k$, but no
jump deformations to algebras on a level higher than $k$. For families,
if one algebra in the family has a jump to an element on level $k$,
then we place the the entire family on at least level $k+1$. Thus, even
though generically, elements of the family $d_{37}(p:q)$ deform only to
members of the same family, there is an element in the family which has
a jump to an element on level 4. For the generic element in a family,
we consider it to be on a higher level than the other elements because
it has jump deformations to the other elements in its family.

\begin{table}[h,t]
\begin{center}
\begin{tabular}{ll}
Level&Codifferentials\\ \hline \\
$1$&1,2,3,4,5,10,11,12,13,14,15,16,17,18,19,20,39,53,54,55,56\\
$2$&6,7,21,22,23,24,25,26,27,28,29,30,40,41,42,43,44,57\\
$3$&8,9,31,32,45,46,47,48,58,59,60,61,62,63,64\\
$4$&33,34,35,36,49,50,66,67,69,70\\
$5$&65,68,$37(p:q)$, $38(p:q)$\\
$6$&$37(0:0)$,$38(0:0),
$ 51,52\\
$7$&71,72\\
\\ \hline
\end{tabular}
\end{center}
\label{Levels}
\caption{The levels of the algebras}
\end{table}

\bibliographystyle{amsplain}
\bibliography{global}

\end{document}